# The Cassels-Tate pairing on polarized abelian varieties

By Bjorn Poonen and Michael Stoll*

**Abstract**

Let $(A, \lambda)$ be a principally polarized abelian variety defined over a global field $k$, and let $\Sha(A)$ be its Shafarevich-Tate group. Let $\Sha(A)_{\mathrm{nd}}$ denote the quotient of $\Sha(A)$ by its maximal divisible subgroup. Cassels and Tate constructed a nondegenerate pairing

$$\Sha(A)_{\mathrm{nd}} \times \Sha(A)_{\mathrm{nd}} \to \mathbf{Q}/\mathbf{Z}.$$

If $A$ is an elliptic curve, then by a result of Cassels the pairing is alternating. But in general it is only antisymmetric.

Using some new but equivalent definitions of the pairing, we derive general criteria deciding whether it is alternating and whether there exists some alternating nondegenerate pairing on $\Sha(A)_{\mathrm{nd}}$. These criteria are expressed in terms of an element $c \in \Sha(A)_{\mathrm{nd}}$ that is canonically associated to the polarization $\lambda$. In the case that $A$ is the Jacobian of some curve, a down-to-earth version of the result allows us to determine effectively whether $\#\Sha(A)$ (if finite) is a square or twice a square. We then apply this to prove that a positive proportion (in some precise sense) of all hyperelliptic curves of even genus $g \geq 2$ over $\mathbf{Q}$ have a Jacobian with nonsquare $\#\Sha$ (if finite). For example, it appears that this density is about 13% for curves of genus 2. The proof makes use of a general result relating global and local densities; this result can be applied in other situations.

## Contents

---

*Much of this research was done while the first author was supported by an NSF Mathematical Sciences Postdoctoral Research Fellowship at Princeton University. This work forms part of the second author's Habilitation in Düsseldorf.

1991 *Mathematics Subject Classification*. 11G10, 11G30, 14H25, 14H40.

## 1. Introduction

The study of the Shafarevich-Tate group $\mathrm{III}(A)$ of an abelian variety $A$ over a global field $k$ is fundamental to the understanding of the arithmetic of $A$. It plays a role analogous to that of the class group in the theory of the multiplicative group over an order in $k$. Cassels [Ca], in one of the first papers devoted to the study of $\mathrm{III}$, proved that in the case where $E$ is an elliptic curve over a number field, there exists a pairing

$$\mathrm{III}(E) \times \mathrm{III}(E) \longrightarrow \mathbf{Q}/\mathbf{Z}$$

that becomes nondegenerate after one divides $\mathrm{III}(E)$ by its maximal divisible subgroup. He proved also that this pairing is alternating; i.e., that $\langle x, x \rangle = 0$ for all $x$. If, as is conjectured, $\mathrm{III}(E)$ is always finite, then this would force its order to be a perfect square. Tate [Ta2] soon generalized Cassels' results by proving that for abelian varieties $A$ and their duals $A^\vee$ in general, there is a pairing

$$\mathrm{III}(A) \times \mathrm{III}(A^\vee) \longrightarrow \mathbf{Q}/\mathbf{Z},$$

that is nondegenerate after division by maximal divisible subgroups. He also proved that if $\mathrm{III}(A)$ is mapped to $\mathrm{III}(A^\vee)$ via a polarization *arising from a k-rational divisor on $A$* then the induced pairing on $\mathrm{III}(A)$ is alternating. But it is known that when dim $A > 1$, a $k$-rational polarization need not come from a $k$-rational divisor on $A$. (See Section 4 for the obstruction.) For principally polarized abelian varieties in general,[1] Flach [Fl] proved that the pairing is antisymmetric, by which we mean $\langle x, y \rangle = -\langle y, x \rangle$ for all $x, y$, which is slightly weaker than the alternating condition.

It seems to have been largely forgotten that the alternating property was never proved in general: in a few places in the literature, one can find the claim that the pairing is always alternating for Jacobians of curves over number fields, for example. In Section 10 we will give explicit examples to show that this is not true, and that $\#\text{III}(J)$ need not be a perfect square even if $J$ is a Jacobian of a curve over $\mathbf{Q}$.[2]

One may ask what properties beyond antisymmetry the pairing has in the general case of a principally polarized abelian variety $(A, \lambda)$ over a global field $k$. For simplicity, let us assume here that $\text{III}(A)$ is finite, so that the pairing is nondegenerate. Flach's result implies that $x \mapsto \langle x, x \rangle$ is a homomorphism $\text{III}(A) \to \mathbf{Q}/\mathbf{Z}$, so by nondegeneracy there exists $c \in \text{III}(A)$ such that $\langle x, x \rangle = \langle x, c \rangle$. Since Flach's result implies $2\langle x, x \rangle = 0$, we also have $2c = 0$ by nondegeneracy. It is then natural to ask, what is this element $c \in \text{III}(A)[2]$ that we have canonically associated to $(A, \lambda)$? An intrinsic definition of $c$ is given in Section 4, and it will be shown[3] that $c$ vanishes (i.e., the pairing is alternating) if and only if the polarization arises from a $k$-rational divisor on $A$. This shows that Tate's and Flach's results are each best possible in a certain sense.

Our paper begins with a summary of most of the notation and terminology that will be needed, and with two definitions of the pairing. (We give two more definitions and prove the compatibility of all four in an appendix.) Sections 4 and 5 give the intrinsic definition of $c$, and show that it has the desired property. (Actually, we work a little more generally: $\lambda$ is not assumed to be principal, and in fact it may be a difference of polarizations.) Section 6 develops some consequences of the existence of $c$; for instance if $A$ is principally polarized and $\text{III}(A)$ is finite, then its order is a square or twice a square according as $\langle c, c \rangle$ equals $0$ or $\frac{1}{2}$ in $\mathbf{Q}/\mathbf{Z}$. We call $A$ *even* in the first case and *odd* in the second case.

---

[1]Actually Flach considers this question in a much more general setting.

[2]This is perhaps especially surprising in light of Urabe's recent results [Ur], which imply, for instance, for the analogous situation of a proper smooth geometrically integral surface $X$ over a finite field $k$ of characteristic $p$, that if the prime-to-$p$ part of $\text{Br}(X)$ is finite, the order of this prime-to-$p$ part *is* a square. (There exist "examples" of nonsquare Brauer groups in the literature, but Urabe explains why they are incorrect.) See Section 11 for more comments on the Brauer group.

[3]The statement of this result needs to be modified slightly if the finiteness of $\text{III}(A)$ is not assumed.

The main goal of Sections 7 and 8 is to translate this into a more down-to-earth criterion for the Jacobian of a genus $g$ curve $X$ over $k$: $\langle c, c \rangle = N/2 \in \mathbf{Q}/\mathbf{Z}$ where $N$ is the number of places $v$ of $k$ for which $X$ has no $k_v$-rational divisor of degree $g - 1$. Section 9 applies this criterion to hyperelliptic curves of even genus $g$ over $\mathbf{Q}$, and shows that a positive proportion $\rho_g$ of these (in a sense to be made precise) have odd Jacobian. It also gives an exact formula for $\rho_g$ in terms of certain local densities, and determines the behavior of $\rho_g$ as $g$ goes to infinity. The result relating the local and global densities is quite general and can be applied to other similar questions. Numerical calculations based on the estimates and formulas obtained give an approximate value of 13% for the density $\rho_2$ of curves of genus 2 over $\mathbf{Q}$ with odd Jacobian.

Section 10 applies the criterion to prove that Jacobians of certain Shimura curves are always even. It gives also a few other examples, including an explicit genus 2 curve over $\mathbf{Q}$ for whose Jacobian we can prove unconditionally that $\langle c, c \rangle = \frac{1}{2}$ and $\mathrm{III} \cong \mathbf{Z}/2\mathbf{Z}$, and another for which $\mathrm{III}$ is finite of square order, but with $\langle \, , \, \rangle$ not alternating on it.

Finally, Section 11 addresses the analogous questions for Brauer groups of surfaces over finite fields, recasting an old question of Tate in new terms.

## 2. Notation

Many of the definitions in this section are standard. The reader is encouraged to skim this section and the next, and to proceed to Section 4.

If $S$ is a set, then $2^S$ denotes its power set.

Suppose that $M$ is an abelian group. For each $n \geq 1$, let $M[n] = \{m \in M : nm = 0\}$. Let $M_{\text{tors}} = \bigcup_{n=1}^{\infty} M[n] = \bigoplus_p M(p)$, where for each prime $p$, $M(p) = \bigcup_{n=1}^{\infty} M[p^n]$ denotes the $p$-primary part of the torsion subgroup of $M$. Let $M_{\text{div}}$ be the maximal divisible subgroup of $M$; i.e., $m$ is in $M_{\text{div}}$ if and only if for all $n \geq 1$ there exists $x \in M$ such that $nx = m$. Denote by $M_{\text{nd}}$ the quotient $M/M_{\text{div}}$. (The subscript nd stands for "nondivisible part.") If

$$\langle \, , \, \rangle : M \times M' \longrightarrow \mathbf{Q}/\mathbf{Z}$$

is a bilinear pairing between two abelian groups, then for any $m \in M$, let $m^{\perp} = \{m' \in M' : \langle m, m' \rangle = 0\}$, and for any subgroup $V \subseteq M$, let $V^{\perp} = \bigcap_{v \in V} v^{\perp}$. When $M = M'$, we say that $\langle \, , \, \rangle$ is *antisymmetric* if $\langle a, b \rangle = -\langle b, a \rangle$ for all $a, b \in M$, and *alternating* if $\langle a, a \rangle = 0$ for all $a \in M$. Note that a bilinear pairing $\langle \, , \, \rangle$ on $M$ is antisymmetric if and only if $m \mapsto \langle m, m \rangle$ is a homomorphism. If a pairing is alternating, then it is antisymmetric, but the converse is guaranteed on $M(p)$ only for odd $p$.

If $k$ is a field, then $\overline{k}$ and $k^{\text{sep}}$ denote algebraic and separable closures, and $G_k = \text{Gal}(\overline{k}/k) = \text{Gal}(k^{\text{sep}}/k)$ denotes its absolute Galois group. If $k$ is a global field, then $M_k$ denotes the set of places of $k$. If moreover $v \in M_k$, then $k_v$ denotes the completion, and $G_v = \text{Gal}(k_v^{\text{sep}}/k_v)$ denotes the absolute Galois group of $k_v$.

Suppose that $G$ is a profinite group acting continuously on an abelian group $M$. We use $C^i(G, M)$ (resp. $Z^i(G, M)$ and $H^i(G, M)$) to denote the group of continuous $i$-cochains (resp. $i$-cocyles and $i$-cohomology classes) with values in the $G$-module $M$. If $k$ is understood, we use $C^i(M)$ as an abbreviation for $C^i(G_k, M)$, and similarly for $Z^i(M)$ and $H^i(M)$. If $\alpha \in C^i(G_k, M)$, then $\alpha_v \in C^i(G_v, M)$ denotes its local restriction.[4] If $v$ is a place of a global field, we use $\text{inv}_v$ to denote the usual monomorphism $H^2(G_v, k_v^{\text{sep}*}) = \text{Br}(k_v) \to \mathbf{Q}/\mathbf{Z}$ (which is an isomorphism if $v$ is nonarchimedean).

Varieties will be assumed to be geometrically integral, smooth, and projective, unless otherwise specified. If $X$ is a variety over $k$, let $k(X)$ denote the function field of $X$. If $K$ is a field extension of $k$, then $X_K$ denotes $X \times_k K$, the same variety with the base extended to $\text{Spec} K$. Let $\text{Div}(X) = H^0(G_k, \text{Div}(X_{k^{\text{sep}}}))$ denote the group of ($k$-rational) Weil divisors on $X$. If $f \in k(X)^*$ or $f \in k(X)^*/k^*$, let $(f) \in \text{Div}(X)$ denote the associated principal divisor. If $D \in \text{Div}(X \times Y)$, let ${}^tD \in \text{Div}(Y \times X)$ denote its transpose. Let $\text{Pic}(X)$ denote the group of invertible sheaves (= line bundles) on $X$, let $\text{Pic}^0(X_{k^{\text{sep}}})$ denote the subgroup of $\text{Pic}(X_{k^{\text{sep}}})$ of invertible sheaves algebraically equivalent to 0, and let $\text{Pic}^0(X) = \text{Pic}(X) \cap \text{Pic}^0(X_{k^{\text{sep}}})$. In older terminology, which we will find convenient to use on occasion, $H^0(\text{Pic}(X_{k^{\text{sep}}}))$ is the group of $k$-rational divisor classes, and $\text{Pic}(X)$ is the subgroup of divisor classes that are actually represented by $k$-rational divisors. Define the Néron-Severi group $\text{NS}(X) = \text{Pic}(X)/\text{Pic}^0(X)$. If $D \in \text{Div}(X)$, let $\mathcal{L}(D) \in \text{Pic}(X)$ be the associated invertible sheaf. Let $\text{Div}^0(X)$ be the subgroup of $\text{Div}(X)$ mapping into $\text{Pic}^0(X)$. If $\lambda \in H^0(\text{NS}(V_{k^{\text{sep}}}))$, let $\text{Pic}^\lambda_{V/k}$ denote the component of the Picard scheme $\text{Pic}_{V/k}$ corresponding to $\lambda$. If $k$ has characteristic $p > 0$, and $\dim V \geq 2$, then these schemes need not be reduced [Mu1, Lect. 27], but in any case the associated reduced scheme $\text{Pic}^\lambda_{V/k,\text{red}}$ is a principal homogeneous space of the Picard variety $\text{Pic}^0_{V/k,\text{red}}$, which is an abelian variety over $k$. Also, $\text{Pic}^\lambda_{V/k,\text{red}}(k^{\text{sep}})$ equals the preimage of $\lambda$ under $\text{Pic}(V_{k^{\text{sep}}}) \to \text{NS}(V_{k^{\text{sep}}})$ with its $G_k$-action.

Let $\mathcal{Z}(X)$ (resp. $\mathcal{Z}^i(X)$) denote the group of 0-cycles on $X$ (resp. the set of 0-cycles of degree $i$ on $X$). For any $i \in \mathbf{Z}$, let $\text{Alb}^i_{X/k}$ denote the degree $i$ component of the Albanese scheme. Let $\mathcal{Y}^0(X)$ denote the kernel of the natural map $\mathcal{Z}^0(X) \to \text{Alb}^0_{X/k}(k)$. Then $\text{Alb}^i_{X/k}$ is a principal homogeneous space of the Albanese variety $\text{Alb}^0_{X/k}$ and its $\overline{k}$-points correspond ($G_k$-equivariantly) to elements of $\mathcal{Z}^i(X_{\overline{k}})$ modulo the action of $\mathcal{Y}^0(X_{\overline{k}})$.

---

[4] We will be slightly sloppy in writing this, because we often will intend the "$M$" in $C^i(G_k, M)$ to be a proper subgroup of the "$M$" in $C^i(G_v, M)$; for instance these two $M$'s may be the $k^{\text{sep}}$-points and $k_v^{\text{sep}}$-points of a group scheme $A$ over $k$, in which case we abbreviate by $C^i(G_k, A)$ and $C^i(G_v, A)$ even though the two $A$'s represent different groups of points.

If $X$ is a curve (geometrically integral, smooth, projective, as usual) and $i \in \mathbf{Z}$, let $\text{Pic}^i(X)$ denote the set of elements of $\text{Pic}(X)$ of degree $i$, and let $\text{Pic}^i_{X/k}$ be the degree $i$ component of the Picard scheme, which is a principal homogeneous space of the Jacobian $\text{Pic}^0_{X/k}$ of $X$. (Since $\dim X = 1$, these are already reduced.) Points on $\text{Pic}^i_{X/k}$ over $\overline{k}$ correspond to divisor classes of degree $i$ on $X_{\overline{k}}$. It will be important to keep in mind that the injection $\text{Pic}^i(X) \to \text{Pic}^i_{X/k}(k)$ is not always surjective. (In other words, $k$-rational divisor classes are not always represented by $k$-rational divisors.)

If $A$ is an abelian variety, then $A^\vee = \text{Pic}^0_{A/k,\text{red}} = \text{Pic}^0_{A/k}$ denotes the dual abelian variety. If $X$ is a principal homogeneous space of $A$, then for each $a \in A(k)$, we let $t_a$ denote the translation-by-$a$ map on $X$, and similarly if $x \in X(k)$, then $t_x$ is the trivialization $A \to X$ mapping $0$ to $x$. If $D \in \text{Div} A$, let $D_x = t_x D \in \text{Div} X$. (Note: if $a \in A(k)$, then $t_a^* D = D_{-a}$.) If $\mathcal{L} \in \text{Pic} X$, then $\phi_\mathcal{L}$ denotes the homomorphism $A \to A^\vee$ mapping $a$ to $t_a^* \mathcal{L} \otimes \mathcal{L}^{-1}$. (There is a natural identification $\text{Pic}^0_{X/k} = A^\vee$.) We may also identify $\phi_\mathcal{L}$ with the class of $\mathcal{L}$ in $\text{NS}(X)$. If $D \in \text{Div} X$, then we define $\phi_D = \phi_{\mathcal{L}(D)}$. A *polarization* on $A$ (defined over $k$) is a homomorphism $A \to A^\vee$ (defined over $k$) which over $k^{\text{sep}}$ equals $\phi_\mathcal{L}$ for some ample $\mathcal{L} \in \text{Pic}(X_{k^{\text{sep}}})$. (One can show that this gives the same concept as the usual definition, in which $\overline{k}$ is used instead of $k^{\text{sep}}$.) A *principal polarization* is a polarization that is an isomorphism.

If $A$ is an abelian variety over a global field $k$, then let $\text{III}(A) = \text{III}(k, A)$ be the Shafarevich-Tate group of $A$ over $k$, whose elements we identify with locally trivial principal homogeneous spaces of $A$ (up to equivalence).

Suppose that $V$ and $W$ are varieties over a field $k$, and that $D$ is a divisor on $V \times W$. If $v \in V(\overline{k})$, let $D(v) \in \text{Div}(W_{\overline{k}})$ be the pullback of $D$ under the map $W \to V \times W$ sending $w$ to $(v, w)$, when this makes sense. For $\mathfrak{a} \in \mathcal{Z}(V_{\overline{k}})$, define $D(\mathfrak{a}) \in \text{Div}(W_{\overline{k}})$ by extending linearly, when this makes sense. If $\mathfrak{a} \in \mathcal{Y}^0(V_{\overline{k}})$ and $D(\mathfrak{a})$ is defined, then $D(\mathfrak{a}) = (f)$ for some function $f$ on $W$. (See the proof of Theorem 10 on p. 171 of [La1, VI, §4].) If in addition $\mathfrak{a}' \in \mathcal{Z}^0(W_{\overline{k}})$, and if $f(\mathfrak{a}')$ is defined, we put $D(\mathfrak{a}, \mathfrak{a}') = f(\mathfrak{a}')$ and say that $D(\mathfrak{a}, \mathfrak{a}')$ is defined. If $\mathfrak{a} \in \mathcal{Y}^0(V_{\overline{k}})$ and $\mathfrak{a}' \in \mathcal{Y}^0(W_{\overline{k}})$, then we may interchange $V$ and $W$ to try to define ${}^t D(\mathfrak{a}', \mathfrak{a})$, and Lang's reciprocity law (p. 171 of [La1] again) states that $D(\mathfrak{a}, \mathfrak{a}')$ and ${}^t D(\mathfrak{a}', \mathfrak{a})$ are defined and equal, provided that $\mathfrak{a} \times \mathfrak{a}'$ and $D$ have disjoint supports.

We let $\mu_\infty$ denote the standard Lebesgue measure on $\mathbf{R}^d$, and let $\mu_p$ denote the Haar measure on $\mathbf{Z}_p^d$ normalized to have total mass 1. For $v = (v_1, v_2, \ldots, v_d) \in \mathbf{Z}^d$, define $|v| := \max_i |v_i|$. If $S \subseteq \mathbf{Z}^d$, then the *density* of $S$ is

defined to be

$$\rho(S) := \lim_{N \to \infty} (2N)^{-d} \sum_{v \in S, |v| \leq N} 1,$$

if the limit exists. Define the *upper density* $\overline{\rho}(S)$ and *lower density* $\underline{\rho}(S)$ similarly, except with the limit replaced by a lim sup or lim inf, respectively.

We will use the notation $\mathbf{A}^d$ and $\mathbf{P}^d$ for $d$-dimensional affine and projective space, respectively.

## 3. Two definitions of the Cassels-Tate pairing

In this section, we present the two definitions of the Cassels-Tate pairing used here. The first is well-known [Mi4]. The second appears to be new, but it was partly inspired by Remark 6.12 on page 100 of [Mi4]. In an appendix we will give two other definitions, and show that all four are compatible.

3.1. *The homogeneous space definition.* Let $A$ be an abelian variety over a global field $k$. Suppose $a \in \text{III}(A)$ and $a' \in \text{III}(A^\vee)$. Let $X$ be the (locally trivial) homogeneous space over $k$ representing $a$. Then $\text{Pic}^0(X_{k^{\text{sep}}})$ is canonically isomorphic as $G_k$-module to $\text{Pic}^0(A_{k^{\text{sep}}}) = A^\vee(k^{\text{sep}})$, so that $a'$ corresponds to an element of $H^1(\text{Pic}^0(X_{k^{\text{sep}}}))$, which we may map to an element $b' \in H^2(k^{\text{sep}}(X)^*/k^{\text{sep}*})$ using the long exact sequence associated with

$$0 \longrightarrow \frac{k^{\text{sep}}(X)^*}{k^{\text{sep}*}} \longrightarrow \text{Div}^0(X_{k^{\text{sep}}}) \longrightarrow \text{Pic}^0(X_{k^{\text{sep}}}) \longrightarrow 0.$$

Since $H^3(k^{\text{sep}*}) = 0$, we may lift $b'$ to an element $f' \in H^2(G_k, k^{\text{sep}}(X)^*)$. Then it turns out that $f'_v \in H^2(G_v, k_v^{\text{sep}}(X)^*)$ is the image of an element $c_v \in H^2(G_v, k_v^{\text{sep}*})$.[5] We define

$$\langle a, a' \rangle = \sum_{v \in M_k} \text{inv}_v(c_v) \in \mathbf{Q}/\mathbf{Z}.$$

See Remark 6.11 of [Mi4] for more details. The obvious advantage of this definition over the others is its simplicity.

If $\lambda : A \to A^\vee$ is a homomorphism, then we define a pairing

$$\langle \, , \, \rangle_\lambda : \text{III}(A) \times \text{III}(A) \longrightarrow \mathbf{Q}/\mathbf{Z}$$

by $\langle a, b \rangle_\lambda = \langle a, \lambda b \rangle$.

3.2. *The Albanese-Picard definition.* Let $V$ be a variety (geometrically integral, smooth, projective, as usual) over a global field $k$. Our goal is to

---

[5]One can compute $c_v$ by *evaluating* $f'_v$ at a point in $X(k_v)$, or more generally at an element of $\mathcal{Z}^1(X_{k_v})$ (provided that one avoids the zeros and poles of $f'_v$).

define a pairing

$$\langle \, , \, \rangle_V : \mathrm{III}(\mathrm{Alb}^0_{V/k}) \times \mathrm{III}(\mathrm{Pic}^0_{V/k,\mathrm{red}}) \longrightarrow \mathbf{Q}/\mathbf{Z}. \tag{1}$$

We will first need a partially-defined $G_k$-equivariant pairing

$$[\, , \,] : \mathcal{Y}^0(V_{k^{\mathrm{sep}}}) \times \mathrm{Div}^0(V_{k^{\mathrm{sep}}}) \longrightarrow k^{\mathrm{sep}*}. \tag{2}$$

Temporarily we work instead with a variety $V$ over a separably closed field $K$. Let $A = \mathrm{Alb}^0_{V/K}$ and $A' = \mathrm{Pic}^0_{V/K,\mathrm{red}}$. Let $\mathfrak{P}$ denote a Poincaré divisor on $A \times A'$. Choose a basepoint $P_0 \in V(K)$ to define a map $\phi : V \to A$. Let $\mathfrak{P}_0 = (\phi, 1)^* \mathfrak{P} \in \mathrm{Div}(V \times A')$. Suppose $y \in \mathcal{Y}^0(V)$ and $D' \in \mathrm{Div}^0(V)$. Choose $z' \in \mathcal{Z}^0(A')$ that sums to $\mathcal{L}(D') \in \mathrm{Pic}^0(V) = A'(K)$. Then $D' - {}^t\mathfrak{P}_0(z')$ is the divisor of some function $f'$ on $V$. Define

$$[y, D'] := f'(y) + \mathfrak{P}_0(y, z')$$

if the terms on the right make sense.

We now show that this pairing is independent of choices made. If we change $z'$, we can change it only by an element $y' \in \mathcal{Y}^0(A')$, and then $[y, D']$ changes by $-{}^t\mathfrak{P}_0(y', y) + \mathfrak{P}_0(y, y') = 0$, by Lang reciprocity. If we change $\mathfrak{P}$ by the divisor of a function $F$ on $A \times A'$, then $[y, D']$ changes by $F(\phi(y) \times z') - F(\phi(y) \times z') = 0$. If $E \in \mathrm{Div}(A')$, and if we change $\mathfrak{P}$ by $\pi_2^* E$ ($\pi_2$ being the projection $A \times A' \to A'$), then $[y, D']$ is again unchanged. By the seesaw principle [Mi3, Th. 5.1], the translate of $\mathfrak{P}_0$ by $(a, 0) \in (A \times A')(K)$ differs from $\mathfrak{P}_0$ by a divisor of the form $(F) + \pi_2^* E$; therefore the definition of $[y, D']$ is independent of the choice of $P_0$. It then follows that if $V$ is a variety over *any* field $k$, then we obtain a $G_k$-equivariant pairing (2).

*Remark.* If $V$ is a curve, then an element of $\mathcal{Y}^0(V_{k^{\mathrm{sep}}})$ is the divisor of a function, and the pairing (2) is simply evaluation of the function at the element of $\mathrm{Div}^0(V_{k^{\mathrm{sep}}})$.

We now return to the definition of (1). It will be built from the two exact sequences

$$0 \longrightarrow \mathcal{Y}^0(V_{k^{\mathrm{sep}}}) \longrightarrow \mathcal{Z}^0(V_{k^{\mathrm{sep}}}) \longrightarrow A(k^{\mathrm{sep}}) \longrightarrow 0, \tag{3}$$

$$0 \longrightarrow \frac{k^{\mathrm{sep}}(V)^*}{k^{\mathrm{sep}*}} \longrightarrow \mathrm{Div}^0(V_{k^{\mathrm{sep}}}) \longrightarrow A'(k^{\mathrm{sep}}) \longrightarrow 0,$$

and the two partially-defined pairings

$$[\, , \,] : \mathcal{Y}^0(V_{k^{\mathrm{sep}}}) \times \mathrm{Div}^0(V_{k^{\mathrm{sep}}}) \longrightarrow k^{\mathrm{sep}*}, \tag{4}$$

$$\mathcal{Z}^0(V_{k^{\mathrm{sep}}}) \times \frac{k^{\mathrm{sep}}(V)^*}{k^{\mathrm{sep}*}} \longrightarrow k^{\mathrm{sep}*},$$

the latter defined by lifting the second argument to a function on $V_{k^{\mathrm{sep}}}$, and "evaluating" on the first argument. We may consider $\frac{k^{\mathrm{sep}}(V)^*}{k^{\mathrm{sep}*}}$ to be a subgroup of $\mathrm{Div}^0(V_{k^{\mathrm{sep}}})$, and then the two pairings agree on $\mathcal{Y}^0(V_{k^{\mathrm{sep}}}) \times \frac{k^{\mathrm{sep}}(V)^*}{k^{\mathrm{sep}*}}$, so there will be no ambiguity if we let $\cup$ denote the cup-product pairing on cochains associated to these. We have also the analogous sequences and pairings for each completion of $k$.

Suppose $a \in \mathrm{III}(A)$ and $a' \in \mathrm{III}(A')$. Choose $\alpha \in Z^1(A(k^{\mathrm{sep}}))$ and $\alpha' \in Z^1(A'(k^{\mathrm{sep}}))$ representing $a$ and $a'$, and lift these to $\mathfrak{a} \in C^1(\mathcal{Z}^0(V_{k^{\mathrm{sep}}}))$ and $\mathfrak{a}' \in C^1(\mathrm{Div}^0(V_{k^{\mathrm{sep}}}))$ so that for all $\sigma, \tau \in G_k$, all $G_k$-conjugates of $\mathfrak{a}_\sigma$ have support disjoint from the support of $\mathfrak{a}'_\tau$. Define

$$\eta := d\mathfrak{a} \cup \mathfrak{a}' - \mathfrak{a} \cup d\mathfrak{a}' \qquad \in C^3(k^{\mathrm{sep}*}).$$

We have $d\eta = 0$, but $H^3(k^{\mathrm{sep}*}) = 0$, so $\eta = d\epsilon$ for some $\epsilon \in C^2(k^{\mathrm{sep}*})$.

Since $a$ is locally trivial, we may for each place $v \in M_k$ choose $\beta_v \in A(k_v^{\mathrm{sep}})$ such that $\alpha_v = d\beta_v$. Choose $\mathfrak{b}_v \in \mathcal{Z}^0(V_{k_v^{\mathrm{sep}}})$ mapping to $\beta_v$ and so that for all $\sigma \in G_k$, all $G_v$-conjugates of $\mathfrak{b}_v$ have support disjoint from the support of $\mathfrak{a}'_\sigma$. Then[6]

$$\gamma_v := (\mathfrak{a}_v - d\mathfrak{b}_v) \cup \mathfrak{a}'_v - \mathfrak{b}_v \cup d\mathfrak{a}'_v - \epsilon_v \qquad \in C^2(G_v, k_v^{\mathrm{sep}*})$$

is a 2-cocycle representing some

$$c_v \in H^2(G_v, k_v^{\mathrm{sep}*}) = \mathrm{Br}(k_v) \xrightarrow{\mathrm{inv}_v} \mathbf{Q}/\mathbf{Z}.$$

Define

$$\langle a, a' \rangle_V = \sum_{v \in M_k} \mathrm{inv}_v(c_v).$$

One checks using $d(x \cup y) = dx \cup y + (-1)^{\deg x} x \cup dy$ [AW, p. 106] that this is well-defined and independent of choices. In proving independence of the choice of $\beta_v$ one uses the local triviality of $a'$, which has not been used so far.

This definition appears to be the most useful one for explicit calculations when $A$ is a Jacobian of a curve of genus greater than 1. This is because the present definition involves only divisors on the curve, instead of $m^2$-torsion or homogeneous spaces of the Jacobian, which are more difficult to deal with computationally.

*Remark.* The reader may have recognized the setup of two exact sequences with two pairings as being the same as that required for the definition of the augmented cup-product (see [Mi4, p. 10]). Here we explain the connection. Suppose that we have exact sequences of $G_k$-modules

$$\begin{array}{ccccccccc} 0 & \longrightarrow & M_1 & \longrightarrow & M_2 & \longrightarrow & M_3 & \longrightarrow & 0, \\ 0 & \longrightarrow & N_1 & \longrightarrow & N_2 & \longrightarrow & N_3 & \longrightarrow & 0, \\ 0 & \longrightarrow & P_1 & \longrightarrow & P_2 & \longrightarrow & P_3 & \longrightarrow & 0 \end{array}$$

---

[6]Since it is only the *difference* of the terms $\mathfrak{a}_v$ and $d\mathfrak{b}_v$ that is in $\mathcal{Y}^0(V_{k^{\mathrm{sep}}})$, it would make no sense to replace $(\mathfrak{a}_v - d\mathfrak{b}_v) \cup \mathfrak{a}'_v$ by $\mathfrak{a}_v \cup \mathfrak{a}'_v - d\mathfrak{b}_v \cup \mathfrak{a}'_v$.

and a $G_k$-equivariant bilinear pairing $M_2 \times N_2 \to P_2$ that maps $M_1 \times N_1$ into $P_1$. The pairing induces a pairing $M_3 \times N_1 \to P_3$. If $a \in \ker\left[H^i(M_1) \to H^i(M_2)\right]$ and $a' \in \ker\left[H^j(N_1) \to H^j(N_2)\right]$, then $a$ comes from some $b \in H^{i-1}(M_3)$, which can be paired with $a'$ to give an element of $H^{i+j-1}(P_3)$. If one changes $b$ by an element $c \in H^{i-1}(M_2)$, then the result is unchanged, since the change can be obtained by pairing $c$ with the (zero) image of $a'$ in $H^j(N_2)$ under the cup-product pairing associated with $M_2 \times N_2 \to P_2$. Thus we have a well-defined pairing[7]

(5) $\quad \ker\left[H^i(M_1) \to H^i(M_2)\right] \times \ker\left[H^j(N_1) \to H^j(N_2)\right] \longrightarrow H^{i+j-1}(P_3).$

If we replace each $M_i$ and $N_i$ by *complexes* with terms in degrees 0 and 1, replace each $P_i$ by a complex with a single term, in degree 1, and replace cohomology by hypercohomology, then we obtain a pairing analogous to (5) defined using the augmented cup-product. One obtains the definition of the Cassels-Tate pairing above by noting that:

1. If $\mathbb{A}_k$ is the adèle ring of $k$,
$$\text{III}(A) = \ker\left[H^1(G_k, A(k^{\text{sep}})) \to H^1(G_k, A(k^{\text{sep}} \otimes_k \mathbb{A}_k))\right]$$
by Shapiro's lemma;

2. The analogous statement holds for $\text{III}(A')$; and

3. If $P_1 = k^{\text{sep}*}$ and $P_2 = (k^{\text{sep}} \otimes_k \mathbb{A}_k)^*$, then the cokernel $P_3$ has $H^2(G_k, P_3) = \mathbf{Q}/\mathbf{Z}$ by class field theory.

## 4. The homogeneous space associated to a polarization

Let $A$ be an abelian variety over a field $k$. To each element of $H^0(\text{NS}(A_{k^{\text{sep}}}))$ we can associate a homogeneous space of $A^\vee$ that measures the obstruction to it arising from a $k$-rational divisor on $A$. From
$$0 \longrightarrow A^\vee(k^{\text{sep}}) \longrightarrow \text{Pic}(A_{k^{\text{sep}}}) \longrightarrow \text{NS}(A_{k^{\text{sep}}}) \longrightarrow 0$$
we obtain the long exact sequence

(6) $\quad 0 \longrightarrow A^\vee(k) \longrightarrow \text{Pic}(A) \longrightarrow H^0(\text{NS}(A_{k^{\text{sep}}}))$
$\longrightarrow H^1(A^\vee(k^{\text{sep}})) \longrightarrow H^1(\text{Pic}(A_{k^{\text{sep}}})).$

(We have $H^0(\text{Pic}(A_{k^{\text{sep}}})) = \text{Pic} A$ because $A(k) \neq \emptyset$.) For $\lambda \in H^0(\text{NS}(A_{k^{\text{sep}}}))$, define $c_\lambda$ to be the image of $\lambda$ in $H^1(A^\vee(k^{\text{sep}}))$. By the proof of Theorem 2 in Section 20 of [Mu3], we have $2\lambda = \phi_{\mathcal{L}}$ where $\mathcal{L} = (1,\lambda)^*\mathcal{P} \in \text{Pic}\,A$ is the pullback of the Poincaré bundle $\mathcal{P}$ on $A \times A^\vee$ by $(1,\lambda) : A \to A \times A^\vee$. Hence $2c_\lambda = 0$.[8]

---

[7] The "diminished cup-product"?!

LEMMA 1. *If $k$ is a local field, then $c_\lambda = 0$ for all $\lambda \in H^0(\text{NS}(A_{k^{\text{sep}}}))$.*

*Proof.* Recall that "Tate[9] local duality" [Ta1] gives a pairing
$$H^0(A(k^{\text{sep}})) \times H^1(A^\vee(k^{\text{sep}})) \longrightarrow H^2(k^{\text{sep}*}) \hookrightarrow \mathbf{Q}/\mathbf{Z}$$
that is perfect, at least after we divide by the connected component on the left in the archimedean case. It can be defined as follows: if $P \in H^0(A(k^{\text{sep}})) = A(k)$ and $z \in H^1(A^\vee(k^{\text{sep}}))$, then we use the long exact sequence associated to
$$0 \longrightarrow \frac{k^{\text{sep}}(A)^*}{k^{\text{sep}*}} \longrightarrow \text{Div}^0(A_{k^{\text{sep}}}) \longrightarrow A^\vee(k^{\text{sep}}) \longrightarrow 0$$
to map $z$ to $H^2\left(\frac{k^{\text{sep}}(A)^*}{k^{\text{sep}*}}\right)$, and "evaluate" the result on a degree zero $k$-rational 0-cycle on $A$ representing $P$ to obtain an element of $H^2(k^{\text{sep}*})$.

Suppose $\lambda \in H^0(\text{NS}(A_{k^{\text{sep}}}))$. By (6), $c_\lambda$ is in the kernel of $H^1(A^\vee) \to H^1(\text{Pic}(A_{k^{\text{sep}}}))$. The long exact sequences associated with

$$\begin{array}{ccccccccc}
0 & \longrightarrow & \dfrac{k^{\text{sep}}(A)^*}{k^{\text{sep}*}} & \longrightarrow & \text{Div}^0(A_{k^{\text{sep}}}) & \longrightarrow & A^\vee(k^{\text{sep}}) & \longrightarrow & 0 \\
& & \| & & \downarrow & & \downarrow & & \\
0 & \longrightarrow & \dfrac{k^{\text{sep}}(A)^*}{k^{\text{sep}*}} & \longrightarrow & \text{Div}(A_{k^{\text{sep}}}) & \longrightarrow & \text{Pic}(A_{k^{\text{sep}}}) & \longrightarrow & 0
\end{array}$$

then show that $c_\lambda$ maps to zero in $H^2\left(\frac{k^{\text{sep}}(A)^*}{k^{\text{sep}*}}\right)$, so that every $P \in A(k)$ pairs with $c_\lambda$ to give 0 in $\mathbf{Q}/\mathbf{Z}$. Hence, by Tate local duality, $c_\lambda = 0$. □

COROLLARY 2. *If $k$ is a global field and $\lambda \in H^0(\text{NS}(A_{k^{\text{sep}}}))$, then $c_\lambda \in \text{III}(A^\vee)[2]$.*

PROPOSITION 3. *If $X$ is a principal homogeneous space of $A$ representing $c \in H^1(A)$, and if $\lambda = \phi_{\mathcal{L}}$ for some $\mathcal{L} \in \text{Pic}\,X$, then $c_\lambda$ is the image of $c$ under the map $H^1(A) \to H^1(A^\vee)$ induced by $\lambda$.*

*Proof.* Pick $P \in X(k^{\text{sep}})$, and pick $D \in \text{Div}(X)$ representing $\mathcal{L}$. Then $D_{-P} \in \text{Div}(A_{k^{\text{sep}}})$ and $\lambda = \phi_{D_{-P}}$, as we see by using $t_P$ to identify $\text{Pic}^0(A_{k^{\text{sep}}})$ with $\operatorname{Pic}^0(X_{k^{\text{sep}}})$. By definition $c_\lambda$ is represented by $\xi \in Z^1(A^\vee)$ where

$$\begin{aligned} \xi_\sigma &:= \text{the class of } {}^\sigma(D_{-P}) - D_{-P} \\ &= \text{the class of } D_{-\sigma P} - D_{-P} \\ &= \text{the class of } (D_{-P})_{-(\sigma P - P)} - D_{-P}, \end{aligned}$$

which by definition represents the image of $c$ under $\phi_{D_{-P}} = \lambda$. □

---

[8]One could also obtain this by using the fact that multiplication by $-1$ on $A$ induces $(-1)^2 = +1$ on $\text{NS}(A_{k^{\text{sep}}})$ and $-1$ on $H^1(A^\vee(k^{\text{sep}}))$.

[9]The archimedean case is related to older results of Witt [Wi]. See p. 221 of [Sc].

COROLLARY 4. *If $(J, \lambda)$ is the canonically polarized Jacobian of a curve $X$, then the element $c_\lambda$ is represented by the principal homogeneous space $\operatorname{Pic}^{g-1}_{X/k} \in H^1(J)$.*

*Proof.* The polarization comes from the theta divisor $\Theta$, which is canonically a $k$-rational divisor on the homogeneous space $\operatorname{Pic}^{g-1}_{X/k}$. □

Combining Corollary 4 with Lemma 1 shows that if $X$ is a curve of genus $g$ over a local field $k_v$, then $X$ has a $k_v$-rational divisor class of degree $g - 1$, a fact originally due to Lichtenbaum [Li].

*Question.* Are all polarizations on an abelian variety $A$ of the form $\phi_{\mathcal{L}}$ for some $\mathcal{L} \in \operatorname{Pic} X$, for some principal homogeneous space $X$ of $A$?

The answer to the question is yes for the canonical polarization on a Jacobian (as mentioned above) or a Prym. (For Pryms, the result is contained in Section 6 of [Mu2], which describes a divisor $\Xi$ on a principal homogeneous space $P^+$ such that $\Xi$ gives rise to the polarization.) One can deduce from Lemma 1 that if $\pi : \tilde{C} \to C$ is an unramified double cover of curves of genus $2g - 1$ and $g$, respectively, with $\tilde{C}$, $C$, and $\pi$ all defined over a local field $k$ of characteristic not 2, then there is a $k$-rational divisor class $\mathcal{D}$ of degree $2g - 2$ on $\tilde{C}$ such that $\pi_*\mathcal{D}$ is the canonical class on $C$.

## 5. The obstruction to being alternating

In this section, we show that if $(A, \lambda)$ is a principally polarized abelian variety over a global field $k$, then $c_\lambda$ (or rather its class in $\text{III}(A^\vee)_{\text{nd}}$) measures the obstruction to $\langle \, , \, \rangle_\lambda$ being alternating. More precisely and more generally, we have the following "pairing theorem":

THEOREM 5. *Suppose that $A$ is an abelian variety over a global field $k$ and $\lambda \in H^0(\operatorname{NS}(A_{k^{\text{sep}}}))$. Then for all $a \in \text{III}(A)$,*

$$\langle a, \lambda a + c_\lambda \rangle = \langle a, \lambda a - c_\lambda \rangle = 0.$$

*Proof.* We will use the homogeneous space definition of $\langle\ ,\ \rangle$. Write $\lambda = \phi_D$ for some $D \in \text{Div}(A_{k^{\text{sep}}})$. Let $X$ be the homogeneous space of $A$ corresponding to $a$.

Now fix $P \in X(k^{\text{sep}})$. For any $\sigma \in G_k$, $\lambda = {}^\sigma\lambda = \phi_{(\sigma D)}$. Thus $\lambda a$ is represented by $\xi \in Z^1(A^\vee)$, where

$$\xi_\sigma := \text{the class of } ({}^\sigma D)_{-(\sigma P - P)} - ({}^\sigma D) \quad \in \text{Pic}^0(A_{k^{\text{sep}}}).$$

Under $t_{(\sigma P)}$, $\xi_\sigma$ corresponds to

$$\xi'_\sigma := \text{the class of } ({}^\sigma D)_P - ({}^\sigma D)_{(\sigma P)} \quad \in \text{Pic}^0(X_{k^{\text{sep}}}).$$

By definition, $c_\lambda$ is represented by $\gamma \in Z^1(A^\vee)$, where

$$\gamma_\sigma := \text{the class of } ({}^\sigma D) - D \quad \in \text{Pic}^0(A_{k^{\text{sep}}}).$$

Under $t_P$, $\gamma_\sigma$ corresponds to

$$\gamma'_\sigma := \text{the class of } ({}^\sigma D)_P - D_P \quad \in \text{Pic}^0(X_{k^{\text{sep}}}).$$

(Recall that the identification $\text{Pic}^0(A_{k^{\text{sep}}}) \cong \text{Pic}^0(X_{k^{\text{sep}}})$ is independent of the trivialization chosen.) Thus $\lambda a - c_\lambda$ is represented by an element of $Z^1(A^\vee)$ corresponding to $\alpha' \in Z^1(\text{Pic}^0(X_{k^{\text{sep}}}))$, where

$$\alpha'_\sigma := \xi'_\sigma - \gamma'_\sigma = \text{the class of } D_P - {}^\sigma(D_P) \quad \in \text{Pic}^0(X_{k^{\text{sep}}}).$$

But $\alpha'$ visibly lifts to an element of $Z^1(\text{Div}^0(X_{k^{\text{sep}}}))$ (that even becomes a coboundary when injected into $Z^1(\text{Div}(X_{k^{\text{sep}}}))$); so the element $b'$ in the definition of Section 3.1 is zero. Hence $\langle a, \lambda a - c_\lambda \rangle = 0$. The other equality follows since $2c_\lambda = 0$. □

*Remark.* One can prove an analogue of Theorem 5 for the Cassels-Tate pairing defined by Flach [Fl] for the Shafarevich-Tate group $\text{III}(M)$ defined by Bloch and Kato for a motive over $\mathbf{Q}$ together with a lattice $M$ in its singular cohomology, under the same assumptions that Flach needs for his antisymmetry result. In fact, the antisymmetry implies the existence of $c$, when the argument in Section 1 is used. (A more desirable solution, however, would be to find an intrinsic definition of $c$ in this context.)

## 6. Consequences of the pairing theorem

In this section, we derive several formal consequences of Theorem 5.

COROLLARY 6. *If $A$ is an abelian variety over a global field $k$ and $\lambda \in H^0(\text{NS}(A_{k^{\text{sep}}}))$, then $\langle\ ,\ \rangle_\lambda$ is antisymmetric.*

*Proof.* The map
$$a \mapsto \langle a, a\rangle_\lambda = \langle a, \lambda a\rangle = -\langle a, c_\lambda\rangle$$
is a homomorphism. □

COROLLARY 7. *Suppose that $A$ is an abelian variety over a global field $k$ and $\lambda \in H^0(\mathrm{NS}(A_{k^{\mathrm{sep}}}))$. Then $\langle\ ,\ \rangle_\lambda$ is alternating if and only if $c_\lambda \in \mathrm{III}(A^\vee)_{\mathrm{div}}$.*

For the rest of this section, we assume that $A$ has a principal polarization $\lambda$, which we fix once and for all. We write $\mathrm{III} = \mathrm{III}(A)$ and let $c = \lambda^{-1}c_\lambda \in \mathrm{III}$. We also drop the subscript $\lambda$ on the pairing $\langle\ ,\ \rangle_\lambda : \mathrm{III} \times \mathrm{III} \to \mathbf{Q}/\mathbf{Z}$, so that Theorem 5 becomes $\langle a, a+c\rangle = \langle a, a-c\rangle = 0$.

Define an endomorphism $a \mapsto a^c$ of the group $\mathrm{III}$ by
$$a^c = \begin{cases} a, & \text{if } \langle a, c\rangle = 0 \\ a+c, & \text{if } \langle a, c\rangle = \tfrac{1}{2}. \end{cases}$$

If $\langle c, c\rangle = 0$, then $a \mapsto a^c$ is an automorphism of order 2. If $\langle c, c\rangle = \tfrac{1}{2}$, then $a \mapsto a^c$ is a projection onto $c^\perp$ with kernel $\{0, c\}$. Define the modified pairing $\langle\ ,\ \rangle^c$ on $\mathrm{III}$ by $\langle a, b\rangle^c = \langle a, b^c\rangle$. By Theorem 5, $\langle\ ,\ \rangle^c$ is alternating.

THEOREM 8. *Let $\bar{c}$ be the image of $c$ in $\mathrm{III}_{\mathrm{nd}}$. The following are equivalent:*

1. *$\langle c, c\rangle = 0$.*

2. *The modified pairing $\langle\ ,\ \rangle^c$ is alternating and nondegenerate on $\mathrm{III}_{\mathrm{nd}}$.*

3. *$\#\mathrm{III}_{\mathrm{nd}}[2]$ is a perfect square.*

4. *$\#\mathrm{III}_{\mathrm{nd}}[n]$ is a perfect square for all $n \geq 1$.*

5. *$\#\mathrm{III}_{\mathrm{nd}}(2)$ is a perfect square.*

6. *$\#\mathrm{III}_{\mathrm{nd}}(p)$ is a perfect square for all primes $p$.*

7. *Either*

   (a) *$\bar{c} = 0$ and $\langle\ ,\ \rangle$ is alternating on $\mathrm{III}_{\mathrm{nd}}$, or*

   (b) *For some $n \geq 1$, there exists a subgroup $V \cong \mathbf{Z}/2^n \times \mathbf{Z}/2^n$ of $\mathrm{III}_{\mathrm{nd}}$ with basis $a, b$ such that $\bar{c} = 2^{n-1}a$, $\langle a, a\rangle = 0$, $\langle a, b\rangle = 2^{-n}$, $\langle b, b\rangle = \tfrac{1}{2}$, and $\mathrm{III}_{\mathrm{nd}} = V \oplus V^\perp$, with $\langle\ ,\ \rangle$ alternating and nondegenerate on $V^\perp$.*

*If these equivalent conditions fail, then the following hold:*

I. $\langle c, c \rangle = \frac{1}{2}$.

II. $\#\mathrm{III}_{\mathrm{nd}}[n]$ is a perfect square for odd $n$, and twice a perfect square for even $n$.

III. $\#\mathrm{III}_{\mathrm{nd}}(p)$ is a perfect square for odd primes $p$, and twice a perfect square for $p = 2$.

IV. $\mathrm{III}_{\mathrm{nd}} = \{0, \bar{c}\} \oplus \bar{c}^{\perp}$, and $\langle \, , \, \rangle$ and $\langle \, , \, \rangle^c$ are alternating and nondegenerate on $\bar{c}^{\perp}$.

*Proof.* Since $2c = 0$, either $\langle c, c \rangle = 0$ or $\langle c, c \rangle = \frac{1}{2}$. Suppose $\langle c, c \rangle = 0$. We already know that $\langle \, , \, \rangle^c$ is alternating, and nondegeneracy follows from the nondegeneracy of $\langle \, , \, \rangle$ on $\mathrm{III}_{\mathrm{nd}}$ and the fact that $a \mapsto a^c$ is an automorphism. Thus 1 implies 2. Clearly 2 implies 3, 4, 5, 6. Also 7 implies 1, so for the equivalence it remains to show that 1 implies 7, and that $\langle c, c \rangle = \frac{1}{2}$ would imply II, III, and IV instead. If $\bar{c} = 0$ we are done by Corollary 7. Otherwise pick the smallest $n \geq 1$ such that $\bar{c} \notin 2^n \mathrm{III}_{\mathrm{nd}}$. Write $\bar{c} = 2^{n-1} a$ for some $a \in \mathrm{III}_{\mathrm{nd}}$. Since $\mathrm{III}_{\mathrm{nd}}[2^n]^{\perp} = 2^n \mathrm{III}_{\mathrm{nd}}$, there exists $b \in \mathrm{III}_{\mathrm{nd}}[2^n]$ such that $\langle b, c \rangle = \frac{1}{2}$. Then $\langle b, b \rangle = \langle b, c \rangle = \langle c, b \rangle = \frac{1}{2}$, and $2^{n-1}\langle a, b \rangle = \langle c, b \rangle = \frac{1}{2}$, and by multiplying $a$ by a suitable element of $(\mathbf{Z}/2^n)^*$, we may assume $\langle a, b \rangle = 2^{-n}$. Let $V$ be the subgroup of $\mathrm{III}_{\mathrm{nd}}$ generated by $a$ and $b$. If $n = 1$, then $a = \bar{c}$, so $\langle a, a \rangle = 0$. If $n > 1$, then $\langle a, a \rangle = \langle a, c \rangle = 2^{n-1}\langle a, a \rangle = 0$, so that $\langle a, a \rangle = 0$ in any case. For $p, q \in \mathbf{Z}$, we have

$$\langle pa + qb, a \rangle = p\langle a, a \rangle + q\langle b, a \rangle = -q/2^n ,$$
$$\langle pa + qb, b \rangle = p\langle a, b \rangle + q\langle b, b \rangle = p/2^n + q/2.$$

If both are zero in $\mathbf{Q}/\mathbf{Z}$, then $q \in 2^n\mathbf{Z}$, and $p \in 2^n\mathbf{Z}$. Hence $V \cong \mathbf{Z}/2^n \times \mathbf{Z}/2^n$, and $\langle \, , \, \rangle$ is nondegenerate on $V$; so $\mathrm{III}_{\mathrm{nd}} = V \oplus V^{\perp}$. This completes the proof that 1 implies 7.

On the other hand, if $\langle c, c \rangle = \frac{1}{2}$, then the nondegeneracy of $\langle \, , \, \rangle$ on $\{0, \bar{c}\}$ implies $\mathrm{III}_{\mathrm{nd}} = \{0, \bar{c}\} \oplus \bar{c}^{\perp}$, and that $\langle \, , \, \rangle$ is alternating and nondegenerate on $\bar{c}^{\perp}$. The rest follows easily. □

COROLLARY 9. *Assume that* III *is finite. Then either*

1. $\langle c, c \rangle = 0$, and there is a finite abelian group $T$ such that $\mathrm{III} \cong T \times T$; in particular, $\#\mathrm{III}$ is a square; or

2. $\langle c, c \rangle = \frac{1}{2}$, and there is a finite abelian group $T$ such that $\mathrm{III} \cong \mathbf{Z}/2\mathbf{Z} \times T \times T$; in particular, $\#\mathrm{III}$ is twice a square.

*If we assume only that* $\Sha_{\mathrm{nd}}$ *is finite, the same conclusions hold with* $\Sha_{\mathrm{nd}}$ *in place of* $\Sha$.

*Remark.* One consequence of Theorem 8 is that if $\langle c, c \rangle = \frac{1}{2}$ then we have one pairing on $\Sha_{\mathrm{nd}}$ that is nondegenerate, and another that is alternating, but we cannot find a pairing that is both!

If $(A, \lambda)$ is a principally polarized abelian variety over a global field $k$, we say that $A$ is *even* if the equivalent conditions 1 through 7 of Theorem 8 hold, and we say that $A$ is *odd* otherwise. Theorem 8 shows that these notions do not depend on the principal polarization $\lambda$ (in cases where there are multiple choices for $\lambda$).

## 7. A formula for Albanese and Picard varieties

Let $V$ be a variety over $k$. In the long exact sequences associated to

$$0 \longrightarrow \frac{k^{\mathrm{sep}}(V)^*}{k^{\mathrm{sep}*}} \longrightarrow \mathrm{Div}(V_{k^{\mathrm{sep}}}) \longrightarrow \mathrm{Pic}(V_{k^{\mathrm{sep}}}) \longrightarrow 0$$

and

$$0 \longrightarrow k^{\mathrm{sep}*} \longrightarrow k^{\mathrm{sep}}(V)^* \longrightarrow \frac{k^{\mathrm{sep}}(V)^*}{k^{\mathrm{sep}*}} \longrightarrow 0,$$

the image of $H^0(\mathrm{Div}(V_{k^{\mathrm{sep}}})) = \mathrm{Div}(V)$ in $H^0(\mathrm{Pic}(V_{k^{\mathrm{sep}}}))$ is $\mathrm{Pic}(V)$, and $H^1(k^{\mathrm{sep}}(V)^*) = 0$, so that upon combining the sequences we obtain an exact sequence

(7) $$0 \longrightarrow \mathrm{Pic}(V) \longrightarrow H^0(\mathrm{Pic}(V_{k^{\mathrm{sep}}})) \xrightarrow{\phi} H^2(k^{\mathrm{sep}*}),$$

and define $\phi$ as shown. (We could also have obtained (7) from the Leray spectral sequence.) Let $\phi_v : H^0(G_v, \mathrm{Pic}(V_{k_v^{\mathrm{sep}}})) \to \mathbf{Q}/\mathbf{Z}$ denote the corresponding homomorphism for the ground field $k_v$, composed with $\mathrm{inv}_v : H^2(G_v, k_v^{\mathrm{sep}*}) \to \mathbf{Q}/\mathbf{Z}$.

PROPOSITION 10. *Let $V$ be a variety over a global field $k$. Let $\langle\ ,\ \rangle_V$ denote the pairing of Section* 3.2. *Suppose that $n$ is an integer for which the homogeneous space $\mathrm{Alb}^n_{V/k}$ of $\mathrm{Alb}^0_{V/k}$ is locally trivial. Suppose $\lambda \in H^0(\mathrm{NS}(V_{k^{\mathrm{sep}}}))$ is such that $\mathrm{Pic}^\lambda_{V/k,\mathrm{red}}$ is locally trivial; choose $\mathcal{L}'_v \in \mathrm{Pic}^\lambda_{V/k,\mathrm{red}}(k_v)$. Then $\phi_v(\mathcal{L}'_v) = 0$ for all but finitely many $v$, and*

$$\langle \mathrm{Alb}^n_{V/k}, \mathrm{Pic}^\lambda_{V/k,\mathrm{red}} \rangle_V = -n \sum_{v \in M_k} \phi_v(\mathcal{L}'_v).$$

*Proof.* It is well-known [CM] that $\phi_v$ is the zero map if $V(k_v) \neq \emptyset$, and the latter holds for all but finitely many $v$ (see Remark 1.6 on p. 249 of [La2]).

Choose $\mathfrak{q} \in \mathcal{Z}^n(V_{k^{\mathrm{sep}}})$ and $D' \in \mathrm{Div}(V_{k^{\mathrm{sep}}})$ mapping into $\mathrm{Pic}^\lambda_{V/k,\mathrm{red}}(k^{\mathrm{sep}})$ so that $^\sigma\mathfrak{q}$ and $^\tau D'$ have disjoint supports for all $\sigma, \tau \in G_k$. We use the notation of Section 3.2. We may take $\mathfrak{a}$ and $\mathfrak{a}'$ so that $\mathfrak{a}_\sigma := {}^\sigma\mathfrak{q} - \mathfrak{q}$ and $\mathfrak{a}'_\sigma := {}^\sigma D' - D'$. Then $d\mathfrak{a} = 0$, $d\mathfrak{a}' = 0$, and $\eta = 0$, so we may take $\varepsilon = 0$. Choose $\beta_v$ and $\mathfrak{b}_v$; then $\gamma_v = (\mathfrak{a}_v - d\mathfrak{b}_v) \cup \mathfrak{a}'_v$.

Choose $L'_v \in \mathrm{Div}(V_{k_v^{\mathrm{sep}}})$ mapping to $\mathcal{L}'_v$, and let $E'_v = D'_v - L'_v$, so that $E'_v \in \mathrm{Div}^0(V_{k_v^{\mathrm{sep}}})$. Then $\mathfrak{a}'_v = dD'_v = dL'_v + dE'_v$. But $(\mathfrak{a}_v - d\mathfrak{b}_v) \cup dE'_v = -d\left[(\mathfrak{a}_v - d\mathfrak{b}_v) \cup E'_v\right]$, so $\gamma_v$ is cohomologous to

$$(\mathfrak{a}_v - d\mathfrak{b}_v) \cup dL'_v = (d(\mathfrak{q}_v - \mathfrak{b}_v)) \cup dL'_v = (d(\mathfrak{q}_v - \mathfrak{b}_v)) \cup f'_v,$$

where $f'_v \in k_v^{\mathrm{sep}}(V)^*$ has divisor $dL'_v$ and we are abusing $\cup$ yet again, this time using it to denote the cup-product pairing associated with the "evaluation pairing"

$$\begin{array}{rcl}\mathcal{Z}(V_{k_v^{\mathrm{sep}}}) \times k_v^{\mathrm{sep}}(V)^* & \longrightarrow & k_v^{\mathrm{sep}*} \\ (\mathfrak{a}\ ,\ f) & \longmapsto & f(\mathfrak{a}).\end{array}$$

(When the first argument is in $\mathcal{Z}^0(V_{k_v^{\mathrm{sep}}})$, this is compatible with the analogue of (4) over $k_v^{\mathrm{sep}}$.)

On the other hand, the image of $\mathcal{L}'_v$ under

$$H^0(G_v, \mathrm{Pic}(V_{k_v^{\mathrm{sep}}})) \longrightarrow H^1\left(G_v, \frac{k_v^{\mathrm{sep}}(V)^*}{k_v^{\mathrm{sep}*}}\right)$$

is represented by $f'_v$ modulo scalars, so $\phi_v(\mathcal{L}'_v)$ is by definition the invariant of $df'_v$ considered as a 2-cocycle taking values in $k_v^{\mathrm{sep}*}$. Since $df'_v$ takes values in $k_v^{\mathrm{sep}*}$ and $\deg(-(\mathfrak{q}_v - \mathfrak{b}_v)) = -(n - 0) = -n$, $-(\mathfrak{q}_v - \mathfrak{b}_v) \cup df'_v \in Z^2(k_v^{\mathrm{sep}*})$ represents $-n\phi_v(\mathcal{L}'_v)$. But $-(\mathfrak{q}_v - \mathfrak{b}_v) \cup df'_v$ differs from the 2-cocycle $(d(\mathfrak{q}_v - \mathfrak{b}_v)) \cup f'_v$ above by the coboundary of $(\mathfrak{q}_v - \mathfrak{b}_v) \cup f'_v$; so they become equal in $H^2(G_v, k_v^{\mathrm{sep}*}) = \mathbf{Q}/\mathbf{Z}$. Summing over $v$ yields the proposition. □

*Remark.* The element $\phi_v(\mathcal{L}'_v) \in \mathbf{Q}/\mathbf{Z}$ may depend on the choice of $\mathcal{L}'_v$, but $n\phi_v(\mathcal{L}'_v)$ does not.

## 8. The criterion for oddness of Jacobians

THEOREM 11. *Suppose that $J$ is the canonically polarized Jacobian of a genus $g$ curve $X$ over a global field $k$. Let $n$ be an integer for which $\mathrm{Pic}^n_{X/k} \in \mathrm{III}(J)$. Then $\langle \mathrm{Pic}^n_{X/k}, \mathrm{Pic}^n_{X/k} \rangle = N/2 \in \mathbf{Q}/\mathbf{Z}$, where $N$ is the number of places $v$ of $k$ for which $\mathrm{Pic}^n(X_{k_v}) = \emptyset$.*

*Proof.* Up to a sign which does not concern us, the pairing $\langle\ ,\ \rangle$ on $\mathrm{III}(J)$ arising from the canonical polarization is the same as the Albanese-Picard pairing for $X$, and $\mathrm{Alb}_{X/k}^n = \mathrm{Pic}_{X/k,\mathrm{red}}^n = \mathrm{Pic}_{X/k}^n$. Let $\mathcal{L}_v$ be a $k_v$-rational divisor class on $X$ of degree $n$. By Proposition 10, it suffices to show that $n\phi_v(\mathcal{L}_v)$ equals $\frac{1}{2}$ or 0 in $\mathbf{Q}/\mathbf{Z}$, according to whether $\mathrm{Pic}^n(X_{k_v}) = \emptyset$ or not.

The *index* $I_v$ (resp. the *period* $P_v$) of $X$ over $k_v$ is the greatest common divisor of the degrees of all $k_v$-rational divisors (resp. $k_v$-rational divisor classes) on $X$. By [Li],[10] $I_v$ equals $P_v$ or $2P_v$, and

$$\begin{aligned}(8) \qquad \phi_v(H^0(G_v, \mathrm{Pic}(X_{k_v^{\mathrm{sep}}}))) &= I_v^{-1}\mathbf{Z}/\mathbf{Z}, \\ \phi_v(H^0(G_v, \mathrm{Pic}^0(X_{k_v^{\mathrm{sep}}}))) &= P_v^{-1}\mathbf{Z}/\mathbf{Z}.\end{aligned}$$

Since $\mathrm{Pic}_{X/k}^n$ is locally trivial, we have $n = r_v P_v$ for some $r_v \in \mathbf{Z}$.

If $\mathrm{Pic}^n(X_{k_v}) \neq \emptyset$ we may take $\mathcal{L}_v \in \mathrm{Pic}^n(X_{k_v})$. By (7) over $k_v$, $\phi(\mathcal{L}_v) = 0$, so $n \cdot \phi_v(\mathcal{L}_v) = 0$, as desired.

If $\mathrm{Pic}^n(X_{k_v}) = \emptyset$ then we must have $I_v = 2P_v$ and $r_v$ odd. Let $\mathcal{D}_v$ denote a $k_v$-rational divisor class of degree $P_v$. We may assume $\mathcal{L}_v = r_v \mathcal{D}_v$. By (8), $\phi_v(\mathcal{D}_v)$ must generate $I_v^{-1}\mathbf{Z}/P_v^{-1}\mathbf{Z}$, so that $P_v\phi_v(\mathcal{D}_v)$ generates $P_v I_v^{-1}\mathbf{Z}/\mathbf{Z} = \frac{1}{2}\mathbf{Z}/\mathbf{Z}$. Hence $P_v\phi_v(\mathcal{D}_v) = \frac{1}{2} \in \mathbf{Q}/\mathbf{Z}$. Then

$$n\phi_v(\mathcal{L}_v) = r_v^2 P_v \phi_v(\mathcal{D}_v) = r_v^2/2 = \tfrac{1}{2} \in \mathbf{Q}/\mathbf{Z}. \qquad \square$$

If $X$ is a curve of genus $g$ over a local field $k_v$, we will say $X$ is *deficient* if $\mathrm{Pic}^{g-1}(X_{k_v}) = \emptyset$, i.e., if $X$ has no $k_v$-rational divisor of degree $g-1$. (Recall from Section 4 that there is always a $k_v$-rational divisor *class* of degree $g-1$.) If $X$ is a curve of genus $g$ over a global field $k$, then a place $v$ of $k$ will be called *deficient* if $X_{k_v}$ is deficient.

COROLLARY 12. *If $J$ is the canonically polarized Jacobian of a genus $g$ curve $X$ over a global field $k$, then the element $c \in \mathrm{III}(J)[2]$ of Section 6 is $\mathrm{Pic}_{X/k}^{g-1}$, and $\langle c,c \rangle = N/2 \in \mathbf{Q}/\mathbf{Z}$, where $N$ is the number of deficient places of $X$.*

*Proof.* Combine Corollary 4 with Theorem 11. $\qquad \square$

*Remarks.* For $g = 1$, Corollaries 7 and 12 let us recover Cassels' theorem [Ca] that $\langle\,,\,\rangle$ is alternating for elliptic curves.

Theorem 11 and Corollary 12 hold even if $X$ has genus 0! In this case $J$ is a point, and so $c = 0$. The condition $\mathrm{Pic}^{-1}(X_{k_v}) = \emptyset$ is equivalent to $X(k_v) = \emptyset$. So we recover the well-known fact that $X(k_v) = \emptyset$ for an even number of places $v$ of $k$.

*Question.* Is there an analogue of Corollary 12 for general principally polarized abelian varieties, or at least an analogue for Pryms?

---

[10] Actually, the proofs in [Li] are for finite extensions of $\mathbf{Q}_p$ only, but the same proofs work for local fields in general, when one uses the duality theorems from [Mi4].

## 9. The density of odd hyperelliptic Jacobians over Q

Fix $g \geq 2$. In this section we are interested in the "probability" that a Jacobian of a random hyperelliptic curve of genus $g$ over $\mathbf{Q}$ has $\langle c, c \rangle = \frac{1}{2}$ (and hence nonsquare order of Ш if Ш is finite).[11] More precisely, we will consider the nonsingular projective models $X$ of curves defined by equations of the form[12]

$$(9) \qquad y^2 = a_{2g+2} x^{2g+2} + a_{2g+1} x^{2g+1} + \cdots + a_1 x + a_0$$

with $a_i \in \mathbf{Z}$.

Let $S_g$ be the set of $a = (a_0, a_1, \ldots, a_{2g+2})$ in $\mathbf{Z}^{2g+3}$ such that (9) defines a hyperelliptic curve of genus $g$ whose Jacobian over $\mathbf{Q}$ is odd. Our goal in this section is to prove that the density $\rho_g := \rho(S_g)$ exists, and to derive an expression for it in terms of certain local densities.

The set of $a \in \mathbf{Z}^{2g+3}$ (resp. in $\mathbf{Z}_p^{2g+3}$ or in $\mathbf{R}^{2g+3}$) for which (9) does *not* define a hyperelliptic curve of genus $g$ is a set of density (resp. measure) zero, because it is the zero locus of the discriminant $\Delta \in \mathbf{Z}[a_0, a_1, \ldots, a_{2g+2}]$. We may hence disregard this set when computing densities or measures.

### 9.1. The archimedean density.
Let $U_{g,\infty}$ be the set of $a \in \mathbf{R}^{2g+3}$ such that the curve over $\mathbf{R}$ defined by (9) has genus $g$ and is deficient. The boundary of $U_{g,\infty}$ is contained in the zero locus of $\Delta$, and hence has measure zero. Let $U^1_{g,\infty} = U_{g,\infty} \cap [-1,1]^{2g+3}$, and let $s_{g,\infty} = 2^{-(2g+3)} \mu_\infty(U^1_{g,\infty})$.

If $g$ is odd, then clearly $s_{g,\infty} = 0$. If $g$ is even, then $a \in U_{g,\infty}$ if and only if the polynomial $\sum_{i=0}^{2g+2} a_i x^i$ is negative on $\mathbf{R}$ and has no double root. Hence for even $g$, $s_{g,\infty} = q_{2g+2}/2$, where $q_n$ denotes the probability that a degree $n$ polynomial with random coefficients drawn independently and uniformly from $[-1, 1]$ has no real root.

PROPOSITION 13. *For all $g \geq 2$, $s_{g,\infty} \leq 1/4$.*

*Proof.* The set $U^1_{g,\infty}$ is contained in the subset of $[-1,1]^{2g+3}$ for which the first and last coordinates are negative. □

PROPOSITION 14. *As $g \to \infty$, $s_{g,\infty} = O(1/\log g)$.*

*Proof.* An old result of Littlewood and Offord [LO1], [LO2] implies that $q_n = O(1/\log n)$ as $n \to \infty$. See also [BS, p. 51]. □

---

[11] We could derive similar results for other number fields, or for other families of curves, but the results would be more awkward to state, so we will restrict attention to hyperelliptic curves over $\mathbf{Q}$.

[12] If $g$ is even, then every hyperelliptic curve has a model of this form. If $g$ is odd, then we are missing the hyperelliptic curves for which the quotient by the hyperelliptic involution is a twisted form of $\mathbf{P}^1$.

Based on some heuristic arguments and Monte Carlo runs, we conjecture that in fact there exists $\alpha > 0$ such that $q_n = n^{-\alpha + o(1)}$ as $n \to \infty$, and $\alpha$ is unchanged if the uniform distribution on $[-1, 1]$ is replaced by any random variable $X$ satisfying

1. $X$ has mean zero,

2. $\text{Prob}(X \neq 0) > 0$, and

3. $X$ belongs to the domain of attraction of the normal law.

Currently there seems to be no such $X$ for which an upper or lower bound of order $n^{-\beta}$ for some $\beta > 0$ is known. On the other hand, there is a huge literature on the *expected number* of real zeros of a random polynomial; see [BS] and [EK] for references. This seems to be an easier question, because of the linearity of expected value (one can sum the expected number of zeros in tiny intervals $[t, t + dt]$).

9.2. *The nonarchimedean densities.* For finite primes $p$, let $U_{g,p}$ be the set of $a \in \mathbf{Z}_p^{2g+3}$ such that the curve over $\mathbf{Q}_p$ defined by (9) has genus $g$ and is deficient. The boundary of $U_{g,p}$ is contained in the zero locus of $\Delta$, and hence has measure zero. In this section we derive bounds on $s_{g,p} := \mu_p(U_{g,p})$.

*Remarks.* Suppose that the curve (9) has genus $g$. Let $\pi : X \to \mathbf{P}^1$ be the projection onto the $x$-line. Then $\pi^* \mathcal{O}(1) \in \text{Pic}^2(X)$. Hence if $g$ is odd, $\text{Pic}^{g-1}(X) \neq \emptyset$, so that $\rho_g = 0$ and $s_{g,p} = 0$ for all $p$.[13] For even $g$, we can only say that the condition $\text{Pic}^{g-1}(X_{\mathbf{Q}_p}) \neq \emptyset$ is equivalent to $\text{Pic}^1(X_{\mathbf{Q}_p}) \neq \emptyset$ (or to $\text{Pic}^m(X_{\mathbf{Q}_p}) \neq \emptyset$ for any other odd integer $m$, for that matter).

The papers [vGY1] and [vGY2] give criteria in special cases for when a hyperelliptic curve of genus $g$ over a $p$-adic field admits a line bundle of odd degree. The first half of the following lemma corresponds to Proposition 3.1 in [vGY1].

LEMMA 15. *Let $X$ be the curve $y^2 = f(x)$ with $f(x) \in \mathbf{Z}[x]$ square free (except possibly for integer squares). Let $\bar{f}$ be the reduction of $f$ modulo an odd prime $p$. Suppose that $\bar{f}$ is* not *of the form $\bar{u}\bar{h}^2$ where $\bar{u} \in \mathbf{F}_p^* \setminus \mathbf{F}_p^{*2}$ and $\bar{h} \in \mathbf{F}_p[x]$ ($\bar{h} = 0$ is allowed). Then $\text{Pic}^1(X_{\mathbf{Q}_p}) \neq \emptyset$. If moreover $p$ is sufficiently large compared to $\deg f$ (or the genus of $X$), then $X(\mathbf{Q}_p) \neq \emptyset$.*

---

[13] Even if we considered hyperelliptic curves for which the quotient by the hyperelliptic involution was a twist of $\mathbf{P}^1$, it would still be true for $g \equiv 1 \pmod{4}$ that a hyperelliptic curve of genus $g$ automatically admitted a line bundle of degree $g - 1$.

*Proof.* Suppose first that $\bar{f}$ is not a constant multiple of a square in $\mathbf{F}_p[x]$ (and so in particular $\bar{f} \neq 0$). Write $\bar{f} = \bar{\ell}\,\bar{j}^2$ with $\bar{\ell}, \bar{j} \in \mathbf{F}_p[x]$ and $\bar{\ell}$ (nonconstant and) square free. Let $\bar{X}$ denote the affine curve $y^2 = \bar{\ell}(x)$ over $\mathbf{F}_p$. By the Weil conjectures, $\#\bar{X}(\mathbf{F}_{p^m}) \to \infty$ as $m \to \infty$, so for any sufficiently large odd $m$, $\bar{X}$ has a point with $x$-coordinate $\bar{\alpha} \in \mathbf{F}_{p^m}$ such that $\bar{\alpha}$ is not a zero of $\bar{f}$. Lift $\bar{\alpha}$ to an element $\alpha$ of the degree $m$ unramified extension $L$ of $\mathbf{Q}_p$. Then $X(L)$ contains a point $P$ with $x$-coordinate $\alpha$. Finally the trace $\mathrm{Tr}_{L/\mathbf{Q}_p} P$ is a $\mathbf{Q}_p$-rational divisor of odd degree $m$. By the remarks preceding the lemma, $\mathrm{Pic}^1(X_{\mathbf{Q}_p}) \neq \emptyset$.

On the other hand, if $\bar{f} = \bar{h}^2$ for a nonzero $\bar{h} \in \mathbf{F}_p[x]$, then for sufficiently large odd $m$, there exists $\bar{\alpha} \in \mathbf{F}_{p^m}$ such that $\bar{h}(\bar{\alpha}) \neq 0$. Lift $\bar{\alpha}$ to an element $\alpha$ of the degree $m$ unramified extension $L$ of $\mathbf{Q}_p$. Then $X(L)$ contains a point with $x$-coordinate $\alpha$, and we conclude as before that $\mathrm{Pic}^1(X_{\mathbf{Q}_p}) \neq \emptyset$.

If $p$ is sufficiently large compared to $\deg f$, then in each case $m = 1$ is sufficient, and so $X(\mathbf{Q}_p) \neq \emptyset$. $\square$

We have a partial converse to Lemma 15.

LEMMA 16. *Let $X$ be a hyperelliptic curve $y^2 = f(x)$ with $f(x) \in \mathbf{Z}[x]$ square free (except possibly for integer squares) and of degree $2g + 2$. Let $p$ be an odd prime. Write $f(x, z)$ for the homogenization of $f(x)$. If $f(x, z) = u h(x, z)^2 + p j(x, z)$ where the reduction of $u \in \mathbf{Z}_p$ is in $\mathbf{F}_p^* \setminus \mathbf{F}_p^{*2}$, and $h(x, z)$ and $j(x, z)$ are homogeneous of degrees $g + 1$ and $2g + 2$, respectively, and if $\bar{h}(x, z)$ and $\bar{j}(x, z)$ (the reductions of $h$ and $j$ mod $p$) have no common factor of odd degree, then $\mathrm{Pic}^1(X_{\mathbf{Q}_p}) = \emptyset$.*

*Proof.* It suffices to show that if $a, b$ are in the ring of integers $\mathcal{O}$ of an odd degree extension $L$ of $\mathbf{Q}_p$, and not both in the maximal ideal $\mathfrak{m}$ of $L$, then $f(a, b)$ is not a square in $L$. Let $v$ denote the (additive) normalized valuation of $L$. If $v(u\,h(a, b)^2) < v(p\,j(a, b))$, then clearly $f(a, b)$ cannot be a square in $L$. On the other hand, if $v(u\,h(a, b)^2) \geq v(p\,j(a, b))$, then $h(a, b)$ must be in $\mathfrak{m}$, and $j(a, b)$ cannot be in $\mathfrak{m}$, since otherwise $(\bar{a} : \bar{b}) \in \mathbf{P}^1(k)$ would be a common zero of $\bar{h}$ and $\bar{j}$, where $k$, the residue field of $L$, is of odd degree over $\mathbf{F}_p$. This would contradict the assumption on $\bar{h}$ and $\bar{j}$. Therefore $v(p\,j(a, b)) = v(p)$ is odd (since $L$ has odd degree over $\mathbf{Q}_p$), which implies $v(u\,h(a, b)^2) > v(p\,j(a, b))$. But then the valuation of $f(a, b)$ is that of $p\,j(a, b)$ and hence odd, so that $f(a, b)$ again cannot be a square. $\square$

*Remarks.* 1. The previous two lemmas are also related to some recent work by Colliot-Thélène and Saito [CTS], and by Bosch and Liu [BL]. Let $k_v$ be a local field with ring of integers $\mathcal{O}_v$ and (finite) residue field $\kappa$. The *index* $I_2$ of a (geometrically integral, smooth, proper) curve $X$ over $k_v$ is the gcd of the degrees of the closed points of $X$. Equivalently, $\deg(\mathrm{Div}(X)) = I_2 \mathbf{Z}$. Let $\mathcal{X} \to \mathrm{Spec}\,\mathcal{O}_v$ denote the minimal proper regular model of $X$. Let $Y_1, \ldots, Y_r$ denote the components of the special fiber of $\mathcal{X}$, and let $d_1, \ldots, d_r$ be their multiplicities. Let $e_i = [\overline{\kappa} \cap \kappa(Y_i) : \kappa]$. Let $I_3 = \gcd\{d_i e_i\}$. Théorème 3.1 of [CTS] implies that $I_2 | I_3$, and Corollary 1.3 of [BL] implies that $I_3 | I_2$. Hence $I_2 = I_3$. This means that $X$ is deficient if and only if $I_3$ does not divide $g - 1$, where $g$ is the genus of $X$. (This is true for arbitrary curves, not only for hyperelliptic ones, where the index is either 1 or 2.)

2. A discussion with J. van Geel yielded the following alternative criterion for deficiency. Let $k$ be $\mathbf{R}$ or a finite extension of $\mathbf{Q}_p$. Let $\mathbf{H}$ be the nontrivial quaternion algebra with center $k$, and let $\mathbf{H}[x]$ be the polynomial ring in which $x$ commutes with scalars (but scalars do not usually commute with each other). Let $X$ be the curve $y^2 = f(x)$ over $k$ with $f$ square free. Then $\mathrm{Pic}^1(X) = \emptyset$ if and only if $f(x)$ is the square of a purely imaginary polynomial in $\mathbf{H}[x]$. ("Purely imaginary" means that the coefficients $c$ satisfy $\iota(c) = -c$ where $\iota$ is the involution of $\mathbf{H}$.) This criterion does not seem to help computationally, however.

PROPOSITION 17. *If $g \geq 2$ is even and $p$ is an odd prime, then*

$$\left(1 - \frac{1}{p}\right)\left(1 - \frac{1}{p^2}\right) \frac{1}{2p^{g+1}} \leq s_{g,p} \leq \left(1 + \frac{1}{p^{g+2}}\right) \frac{1}{2p^{g+1}}.$$

*Proof.* An easy counting argument shows that exactly $\frac{1}{2}(p^{g+2} - 1) + 1$ of the $p^{2g+3}$ polynomials $\bar{f} \in \mathbf{F}_p[x]$ of degree at most $2g + 2$ are of the form $\bar{u}\bar{h}^2$ with $\bar{u} \in \mathbf{F}_p^* \setminus \mathbf{F}_p^{*2}$ and $\bar{h} \in \mathbf{F}_p[x]$. Together with Lemma 15, this proves the upper bound.

Similarly, the probability that the homogenization $\bar{f}(x, z)$ is of the form $\bar{u}\bar{h}^2$ with $\bar{u} \in \mathbf{F}_p^* \setminus \mathbf{F}_p^{*2}$ and $\bar{h}(x, z) \in \mathbf{F}_p[x, z]$ homogeneous and nonzero is $(p^{-(g+1)} - p^{-(2g+3)})/2$. In this situation, lift $\bar{u}$ to $u \in \{1, 2, \ldots, p - 1\}$, lift $\bar{h}$ to $h$ with coefficients in $\{0, 1, \ldots, p - 1\}$, and write $f = uh^2 + pj$, with $j \in \mathbf{Z}[x, z]$ homogeneous of degree $2g + 2$. Assume that $\bar{u}$ was chosen at random uniformly from the nonquadratic residues, and that the sign of $\bar{h}$ was also chosen at random with equal probabilities; the choices do not affect whether $\bar{h}$ and $\bar{j}$ have a common factor. Then the distribution of $\bar{j}$ is uniform and independent of $\bar{h}$. In the following paragraph we show that the probability for two random homogeneous polynomials in two variables to have no common factor is $(1 - p^{-1})(1 - p^{-2})$. Since $\bar{h}$ is already known to be nonzero, the probability that $\bar{h}$ and $\bar{j}$ have no common factor is $(1 - p^{-1})(1 - p^{-2})/(1 - p^{-(g+2)})$. Together with Lemma 16, this proves the lower bound.

Let us now prove that the probability for two (chosen uniformly at random) homogeneous polynomials in $\mathbf{F}_p[x, z]$ of degrees $m, n > 0$ to have no common factor is $(1 - p^{-1})(1 - p^{-2})$. Let $a_n = p^{n+1} - 1$ be the number of nonzero homogeneous polynomials of degree $n$, and let $b_{m,n}$ denote the number of ordered pairs of nonzero homogeneous polynomials of degrees $m$ and $n$ without common factor (if $n$ is negative, we set $a_n = 0$, and similarly for $b_{m,n}$). Then we have the relations

$$a_n - (p+1)\,a_{n-1} + p\,a_{n-2} = (p-1)\,\delta_n\,,$$

$$\sum_j b_{m-j, n-j}\, a_j = (p-1)\,a_m\,a_n\,.$$

In the first equation, $\delta_n = 1$ if $n = 0$ and $\delta_n = 0$ otherwise. In the second equation, the $j^{\text{th}}$ term on the left counts the pairs with a gcd of degree $j$; each pair is counted $(p-1)$ times, since the gcd is determined only mod $\mathbf{F}_p^*$. Now, if $m, n > 0$,

$$\begin{aligned}
(p-1)(p^2-1)p^{m+n-1} &= a_m\,a_n - (p+1)\,a_{m-1}\,a_{n-1} + p\,a_{m-2}\,a_{n-2} \\
&= \sum_j b_{m-j, n-j}\,(a_j - (p+1)\,a_{j-1} + p\,a_{j-2})/(p-1) \\
&= b_{m,n}\,.
\end{aligned}$$

Since the polynomials have to be nonzero in order to have no common factor, the probability we want is $b_{m,n}/p^{m+n+2} = (1 - p^{-1})(1 - p^{-2})$. □

*Remarks.* Proposition 17 implies that for fixed even $g \geq 2$, $s_{g,p} \sim \frac{1}{2}p^{-(g+1)}$ as $p \to \infty$.

A more careful analysis shows that for odd $p$,

$$\tfrac{1}{2}p^{-3}\left(1 - p^{-1} + \tfrac{1}{2}p^{-2} - 2p^{-3} + p^{-4} + p^{-5} - \tfrac{1}{2}p^{-6}\right)$$
$$\leq s_{2,p} \leq \tfrac{1}{2}p^{-3}\left(1 - p^{-1} + p^{-2} + 3p^{-4} - p^{-5}\right).$$

For practical purposes (dealing with a given curve), the following criterion might be useful.

LEMMA 18. *Let $X$ be the curve $y^2 = f(x)$ with $f(x) \in \mathbf{Z}[x]$ square free (except possibly for integer squares) of degree $2g + 2$ with $g$ even, and let $p$ be an odd prime. Then $p$ can be deficient for $X$ only if the discriminant of $f$ is divisible by $p^{g+1}$.*

*Proof.* This follows from Lemma 15 in the following way. Let $F(x, z) = z^{2g+2} f(x/z)$ be the homogenization of $f$, and let $\bar{F}$ denote the reduction of $F$ mod $p$. If $\bar{F} = 0$, then $\operatorname{disc}(f)$ is divisible by $p^{4g+2}$, hence *a fortiori* by $p^{g+1}$.

Now assume that $\bar{F} \neq 0$. Let
$$\bar{F} = \bar{u} \cdot h_1^{e_1} \cdots h_m^{e_m}$$
be its factorization with $h_j$ of degree $d_j$ and $\bar{u} \in \mathbf{F}_p^*$. To this corresponds a factorization over $\mathbf{Z}_p$
$$F = u \cdot H_1 \cdots H_m$$
with $H_j$ of degree $e_j d_j$ and reducing to $h_j^{e_j}$. The valuation of the discriminant of $H_j$ must be at least
$$2 \frac{1}{e_j} d_j \binom{e_j}{2} = d_j (e_j - 1),$$
and hence we get
$$v_p(\operatorname{disc}(f)) = \sum_j v_p(\operatorname{disc}(H_j)) \geq \sum_j d_j (e_j - 1) = \deg(f) - \sum_j d_j.$$

If $\bar{F}$ is (a constant times) a square, then the sum of the $d_j$ is at most $\frac{1}{2} \deg(F)$, whence $p^{g+1} \mid \operatorname{disc}(F) = \operatorname{disc}(f)$. □

PROPOSITION 19. *For all $g \geq 2$, $s_{g,2} \leq 875/1944 \approx 0.4501$ and $s_{g,2} \leq \exp(-0.06 g / \log g)$.*

*Proof.* Suppose $P_1, P_2, \ldots, P_r$ are distinct closed points of $\mathbf{P}^1_{\mathbf{F}_2}$ such that the residue degree $d_i$ of $P_i$ is odd, and such that $\sum_{i=1}^r d_i \leq 2g + 3$. Let $\eta_j$ denote the normalized Haar measure of the nonsquares in the ring of integers of the unramified extension of $\mathbf{Q}_2$ of degree $j$. We claim that $s_{g,2} \leq \prod_{i=1}^r \eta_{d_i}$.

For each $i$, let $\mathcal{O}_i$ be the ring of integers in the $d_i^{\text{th}}$ degree unramified extension of $\mathbf{Q}_2$, and choose a point $Q_i \in \mathbf{A}^2(\mathcal{O}_i)$ with relatively prime coordinates such that the projection of $Q_i$ onto $\mathbf{P}^1(\mathcal{O}_i)$ reduces to a point in $\mathbf{P}^1(\overline{\mathbf{F}}_2)$ above $P_i$. For $a \in \mathbf{Z}_2^{2g+3}$, let $f(x, z)$ denote $\sum_{j=0}^{2g+2} a_j x^j z^{2g+2-j}$. The set $U_{g,2}$ is contained in the set $U'$ of $a \in \mathbf{Z}_2^{2g+3}$ for which $f(Q_i)$ is a nonsquare in $\mathcal{O}_i$ for all $i$. Applying Lagrange interpolation to the $Q_i$ and their $\mathbf{Q}_2$-conjugates shows that the $\mathbf{Z}_2$-module homomorphism
$$\begin{aligned} \mathbf{Z}_2^{2g+3} &\longrightarrow \prod_{i=1}^r \mathcal{O}_i \\ a &\longmapsto (f(Q_1), f(Q_2), \ldots, f(Q_r)) \end{aligned}$$
is surjective. (Remember that $2g + 3 \geq \sum_{i=1}^r d_i$.) Hence
$$s_{g,2} = \mu_2(U_{g,2}) \leq \mu_2(U') = \prod_{i=1}^r \eta_{d_i},$$
as claimed.

A short calculation shows that
$$\eta_j = \frac{3 \cdot 2^{j-2} + 1}{2^j + 1} \leq 5/6$$
for all $j \geq 1$. By taking $P_1$, $P_2$, $P_3$ to be the points in $\mathbf{P}^1(\mathbf{F}_2)$, and $P_4$ to be a point of degree 3, we obtain $s_{g,2} \leq \eta_1^3 \eta_3 = (5/6)^3 (7/9) = 875/1944$.

If instead we let $m$ be the largest odd integer with $2^m \leq 2g + 2$, and take $P_1, P_2, \ldots, P_r$ to be all closed points on $\mathbf{P}^1_{\mathbf{F}_2}$ of degree dividing $m$, then $\sum_{i=1}^r d_i = \#\mathbf{P}^1(\mathbf{F}_{2^m}) = 2^m + 1 \leq 2g + 3$, and
$$r \geq \frac{2^m}{m} \geq \frac{(2g+2)/4}{\log_2[(2g+2)/4]} \geq \frac{g}{2 \log_2 g},$$
so that
$$s_{g,2} \leq (5/6)^r \leq \exp(\alpha g / \log g)$$
where $\alpha = \frac{1}{2} \log(5/6) \log 2 \approx -0.0631878$. □

*Remarks.* 1. It is a fact that $s_{g,p} \in \mathbf{Q}$ for all even $g \geq 2$ and finite $p$. Indeed, Theorem 7.4 in [De][14] implies that any subset of $\mathbf{Z}_p^{2g+3}$ that is definable (in the sense of (6.1) of [De]) has rational Haar measure. To see that $U_{g,p}$ is definable, observe that the Riemann-Roch theorem implies that the existence of a $\mathbf{Q}_p$-rational divisor of degree 1 is equivalent to the existence of a point defined over some extension of $\mathbf{Q}_p$ of odd degree $\leq g + 1$, and the latter is clearly a definable property.

We expect that for fixed even $g \geq 2$, there even exists $N_g \geq 1$ such that for any residue class modulo $N_g$, the values of $s_{g,p}$ for sufficiently large primes in that class are given by a rational function of $p$. Some general uniformity results in $p$ for definable sets can be found in [Mac2] (Thm. 22 and its corollaries), but obtaining rational functions of $p$ will require use of properties of the sets $U_{g,p}$ beyond their definability.

2. In [St], we will describe a practical algorithm deciding whether a place is deficient for a given hyperelliptic curve (of even genus) over a number field. This algorithm then makes it possible to find out whether the curve in question has even or odd Jacobian, which (in the odd case) will improve the bound on the rank given by the 2-descent procedure (see [PS]) and (under the standard finiteness conjecture regarding Ш) will pin down its parity. For curves of genus 2 over $\mathbf{Q}$, the 2-descent, together with the deficiency-checking algorithm, has been implemented as a computer program by the second author. For details, see again [St].

---

[14] It would be more accurate to say "Theorem 7.4 and its proof," because the theorem states that a certain function $Z_S(s)$ is a rational function of $p^{-s}$, whereas we need to know also that the coefficients of the rational function are themselves in $\mathbf{Q}$. Denef's proof uses the quantifier elimination of [Mac1].

9.3. *The passage from local to global.* We formalize the method for obtaining density results in the following lemma, because the method undoubtedly has other applications. (See the remark after Lemma 21 and the remarks after Theorem 23, for instance.)

LEMMA 20. *Suppose that $U_\infty$ is a subset of $\mathbf{R}^d$ such that $\mathbf{R}^+ \cdot U_\infty = U_\infty$ and $\mu_\infty(\partial U_\infty) = 0$. Let $U_\infty^1 = U_\infty \cap [-1,1]^d$, and let $s_\infty = 2^{-d}\mu_\infty(U_\infty^1)$.[15] Suppose that for each finite prime $p$, $U_p$ is a subset of $\mathbf{Z}_p^d$ such that $\mu_p(\partial U_p) = 0$. Let $s_p = \mu_p(U_p)$. Finally, suppose that*

(10)
$$\lim_{M \to \infty} \overline{\rho}\left(\left\{a \in \mathbf{Z}^d : a \in U_p \text{ for some finite prime } p \text{ greater than } M\right\}\right) = 0.$$

*Define a map $P: \mathbf{Z}^d \to 2^{M_\mathbf{Q}}$ as follows: if $a \in \mathbf{Z}^d$, let $P(a)$ be the set of places $p$ such that $a \in U_p$. Then*

1. *$\sum_p s_p$ converges.*

2. *For $\mathfrak{S} \subseteq 2^{M_\mathbf{Q}}$, $\nu(\mathfrak{S}) := \rho(P^{-1}(\mathfrak{S}))$ exists, and $\nu$ defines a measure on $2^{M_\mathbf{Q}}$.*

3. *The measure $\nu$ is concentrated at the finite subsets of $M_\mathbf{Q}$: for each finite subset $S$ of $M_\mathbf{Q}$,*

(11)
$$\nu(\{S\}) = \prod_{p \in S} s_p \prod_{p \notin S}(1 - s_p),$$

*and if $\mathfrak{S} \subset 2^{M_\mathbf{Q}}$ consists of infinite subsets of $M_\mathbf{Q}$, then $\nu(\mathfrak{S}) = 0$.*

*Proof.* First we handle the case in which there exists $M$ such that $U_p = \emptyset$ for all finite $p > M$. For convenience we say that a subset of $\mathbf{Z}_p$ is an *open interval* if it has the form $\{x \in \mathbf{Z}_p : |x - a|_p \leq b\}$ for some $a \in \mathbf{Z}_p$ and $b \in \mathbf{R}$, and a subset of $[-1, 1]$ is an *open interval* if it has the form $U \cap [-1, 1]$ for some open interval $U$ in $\mathbf{R}$. By an *open box* in $[-1, 1]^d \times \prod_{p \leq M} \mathbf{Z}_p^d$ we mean a Cartesian product of open intervals.

Let $B_N = \{a \in \mathbf{Z}^d : |a| \leq N\}$. Embed $B_N$ in $[-1, 1]^d$ by dividing by $N$. We claim that the fraction of $B_N$ that maps in $[-1, 1]^d \times \prod_{p \leq M} \mathbf{Z}_p^d$ into an open box $I = I_\infty \times \prod_{p \leq M} I_p$ approaches $2^{-d}\mu(I)$ as $N \to \infty$, where $\mu$ denotes the product measure of the $\mu_p$. The subset of $B_N$ satisfying the $p$-adic conditions is the subset contained in a finite disjoint union of translates of a sublattice determined by congruence conditions; the number of translates as a fraction of the index of the sublattice equals $\prod_{p \leq M} \mu_p(I_p)$. For each of these translates, the fraction of elements satisfying the box condition at $\infty$ equals

$$\frac{\mu_\infty(I_\infty)}{\mu_\infty([-1,1]^d)} + O\left(\frac{1}{N}\right),$$

which tends to $2^{-d}\mu_\infty(I_\infty)$ as $N \to \infty$, proving our claim.

---

[15]Since $U_\infty^1$ is the union of the open set $(U_\infty^1)^0$ (its interior) and a subset of a measure zero set, $U_\infty^1$ is automatically measurable.

Let $V_p$ be the complement of $U_p$ in $\mathbf{Z}_p^d$ (or of $U_\infty^1$ in $[-1,1]^d$ if $p = \infty$). If for each $p \leq M$ and $p = \infty$ we choose either $U_p$ ($U_\infty^1$ if $p = \infty$) or $V_p$, and let $P$ be the product, then by compactness we can cover $\bar{P}$ (the closure of $P$) by a finite number of open boxes the sum of whose measures is arbitrarily close to the measure of $\bar{P}$, which equals the product of the measures of the $U_p$ or $V_p$, since each has boundary of measure zero. By the previous paragraph it follows that for any $\varepsilon > 0$, the fraction of $B_N$ that maps into $P$ can be bounded above by $2^{-d}\mu(P) + \varepsilon$ for sufficiently large $N$. On the other hand, the fraction of $B_N$ that maps *outside* $P$ can also be bounded above, since the complement of $P$ is a finite disjoint union of other such $P$'s, and the bound obtained is of the form

$$1 - 2^{-d}\mu(P) + (2^L - 1)\varepsilon$$

for sufficiently large $N$, where $L$ equals one plus the number of primes up to $M$. Thus the fraction of $B_N$ that maps into $P$ tends to $2^{-d}\mu(P)$ as $N \to \infty$. This completes the proof in the case where $U_p = \emptyset$ for $p > M$.

We return to the general case: if $M < M'$ then the density of the set of $a \in \mathbf{Z}^d$ for which $a \notin U_p$ for all primes $p$ with $M < p < M'$ is $\prod_{M < p < M'}(1 - s_p)$, and the condition (10) implies that this tends to 1 as $M \to \infty$, uniformly with respect to $M'$. In other words, $\prod_p(1-s_p)$ converges. Thus $\sum_p s_p$ converges. This implies the convergence of the infinite product in (11).

Fix $\mathfrak{S} \subseteq 2^{M_\mathbf{Q}}$. The condition (10) implies that the actual upper and lower densities of $P^{-1}(\mathfrak{S})$ are approximately what we would have obtained for the density if $U_p$ had been replaced by the empty set for finite $p > M$, and that the error in this approximation tends to zero as $M \to \infty$. The other results are now clear. □

To verify that Lemma 20 applies in our situation, we must check (10). This we will do in Lemma 22 below, with the aid of a result of Ekedahl.

LEMMA 21. *Suppose $f$ and $g$ are relatively prime polynomials in $\mathbf{Z}[x_1, x_2, \ldots, x_d]$. Let $S_M(f,g)$ be the set of $a \in \mathbf{Z}^d$ for which there exists a finite prime $p > M$ dividing both $f(a)$ and $g(a)$. Then $\lim_{M \to \infty} \overline{\rho}(S_M(f,g)) = 0$.*

*Proof.* Apply Theorem 1.2 of [Ek] to the subscheme of $\mathbf{A}_\mathbf{Z}^d$ defined by the equations $f = g = 0$. □

*Remark.* Once one has Lemma 21, it is easy to apply Lemma 20 to obtain a formula for the density of $a \in \mathbf{Z}^d$ such that $f(a)$ and $g(a)$ are relatively prime, in terms of the number of solutions to $f(a) \equiv g(a) \equiv 0$ in $\mathbf{F}_p^d$ for each $p$. The same can of course be done for $\{a \in \mathbf{Z}^d : \gcd(f_1(a), \ldots, f_n(a)) = 1\}$, provided that the polynomials $f_i \in \mathbf{Z}[x_1, \ldots, x_d]$ define a subvariety of codimension at least 2 in $\mathbf{A}_\mathbf{C}^d$. Results such as these follow also from [Ek].

LEMMA 22. *Fix $g \geq 1$. Let $R_M$ (resp. $R'_M$) be the set of $a = (a_0, a_1, \ldots, a_{2g+2}) \in \mathbf{Z}^{2g+3}$ for which (9) is a curve $X$ of genus $g$ that fails to admit a $\mathbf{Q}_p$-rational point (resp. a $\mathbf{Q}_p$-rational divisor of degree 1 — for even $g$, this is equivalent to being deficient) at some finite prime $p$ greater than $M$. Then $\lim_{M \to \infty} \overline{\rho}(R_M) = \lim_{M \to \infty} \overline{\rho}(R'_M) = 0$.*

*Proof.* The space $\mathrm{Pol}_n$ of binary homogeneous polynomials $\sum_{i=0}^n a_i x^i z^{n-i}$ over $\mathbf{C}$ may be identified with $\mathbf{A}^{n+1}$. If $g \geq 1$, then the Zariski closure[16] $V$ of the image of the squaring map $\mathrm{Pol}_{g+1} \to \mathrm{Pol}_{2g+2}$ is of codimension at least 2 in $\mathrm{Pol}_{2g+2} = \mathbf{A}^{2g+3}$, so that we can find two relatively prime polynomials $f, g \in \mathbf{Z}[a_0, \ldots, a_{2g+2}]$ that vanish on $V$. For all but finitely many primes $p$, it is true that if $a \in \mathbf{Z}^{2g+3}$ and $\sum_{i=0}^{2g+2} a_i x^i \bmod p$ is a square in $\overline{\mathbf{F}}_p[x]$ then $p$ divides $f(a)$ and $g(a)$. Combining this with Lemma 15 shows that $R_M \subseteq S_M(f, g)$ for sufficiently large $M$. By Lemma 21, $\lim_{M \to \infty} \overline{\rho}(S_M(f, g)) = 0$, so $\lim_{M \to \infty} \overline{\rho}(R_M) = 0$ also. Since $R'_M \subseteq R_M$, $\lim_{M \to \infty} \overline{\rho}(R'_M) = 0$ as well. □

9.4. *The global density.*

THEOREM 23. *For each $g \geq 1$, the density $\rho_g = \rho(S_g)$ exists, and*

(12) $$1 - 2\rho_g = \prod_{p \in M_\mathbf{Q}} (1 - 2s_{g,p}).$$

*If $g$ is odd, then $\rho_g = 0$ and $s_{g,p} = 0$ for all $p$. If $g$ is even, then $0 < \rho_g < 1$, and $0 < s_{g,p} < 1$ for all $p$.*

*Proof.* Let $\mathfrak{S}_{\mathrm{odd}} \subseteq 2^{M_\mathbf{Q}}$ (resp. $\mathfrak{S}_{\mathrm{even}}$) be the set of all finite subsets of $M_\mathbf{Q}$ of odd (resp. even) order. Lemma 22 gives us the condition (10) needed for the application of Lemma 20. Thus $\nu(\mathfrak{S}_{\mathrm{odd}})$ exists. By Corollary 12, $\rho_g = \nu(\mathfrak{S}_{\mathrm{odd}})$. Also by Lemma 20, $\nu(\mathfrak{S}_{\mathrm{even}}) = 1 - \rho_g$, and the quantity $1 - 2\rho_g = \nu(\mathfrak{S}_{\mathrm{even}}) - \nu(\mathfrak{S}_{\mathrm{odd}})$ can be obtained by substituting $x_p = -1$ for all places $p$ in the generating function $\prod_p (1 - s_{g,p} + s_{g,p} x_p)$. The final vanishing results were noted already in Section 9.2. Finally suppose $g$ is even. It is enough, by openness, to exhibit for each $p$ two nonsingular curves over $\mathbf{Q}_p$ in the correct form, one deficient and one not. The curve $y^2 = x^{2g+2} + 1$ is nondeficient at all $p$, so it remains to construct deficient curves. For $p = \infty$ we may take $y^2 = -(x^{2g+2} + 1)$. For odd finite $p$, Lemma 16 gives deficient curves. For $p = 2$, a similar argument shows that $y^2 = 2x^{2g+2} + 3$ is deficient. □

---

[16] In fact, it is easy to show that the image of the squaring map is already Zariski closed, but we do not need this.

*Remarks.* 1. One could also think of (12) as arising from a Fourier transform of an infinite convolution on the group $\mathbf{Z}/2\mathbf{Z}$.

2. We can also deduce from Lemmas 20 and 22 that for any even $g \geq 2$, and any finite subset $S \subseteq M_{\mathbf{Q}}$, there exist infinitely many hyperelliptic curves $X$ over $\mathbf{Q}$ of genus $g$ for which the set of deficient places is exactly $S$. The same can be proved for any number field $k$, provided that $S$ contains no complex places.

Similarly, Lemmas 20 and 22 (generalized to higher number fields) can be used to show that for any $g \geq 1$, any number field $k$, and any finite subset $S \subseteq M_k$ containing no complex places, there exist infinitely many (hyperelliptic[17]) curves $X$ over $k$ of genus $g$ for which $\{\, v \in M_k \mid X(k_v) = \emptyset \,\} = S$.

Another corollary is that for every $g \geq 2$ and every number field $k$, there exists an *odd* principally polarized abelian variety $A$ over $k$ of dimension $g$. To construct $A$, start with a genus 2 curve over $k$ whose Jacobian $J$ is odd, and let $A = J \times E^{g-2}$ for an elliptic curve $E$ over $k$, with the product polarization.

3. It appears that $\rho_2 \approx 0.13$. This value was obtained in the following way. A numerical Monte Carlo integration (with $10^{11}$ tries) gives the value $s_{2,\infty} \approx 0.0983$. In the same way (again with $10^{11}$ tries, using a computer program that checks for points of odd degree over $\mathbf{Q}_2$) we obtain the estimate $s_{2,2} \approx 0.02377$. (Note that this is much smaller than the bound in Proposition 19.) The estimate for $\rho_2$ then follows from Theorem 23 and the bounds in the remarks after Proposition 17.

THEOREM 24. *As $g \to \infty$, $\rho_g = O(1/\log g)$.*

*Proof.* This follows from Propositions 14, 17, and 19, and Theorem 23. □

Our conjecture that $q_n = n^{-\alpha + o(1)}$ as $n \to \infty$, would imply $\rho_g = g^{-\alpha + o(1)}$ for even $g$ tending to infinity.

## 10. Examples of Shafarevich-Tate groups of Jacobians

In this section we apply Theorem 11 to deduce facts about Shafarevich-Tate groups of Jacobians.

---

[17]For $g = 1$, the curves will be of the form $y^2 = f(x)$ with $\deg f = 4$.

### 10.1. Jacobians of Shimura curves.

THEOREM 25. *Let $k$ be a number field, and let $B$ be an indefinite quaternion algebra over $\mathbf{Q}$ of discriminant $\mathrm{Disc}\, B > 1$. Let $X$ be the corresponding Shimura curve, considered over $k$, and let $J$ be its Jacobian. If $d \in \mathbf{Z}$ and $\mathrm{Pic}^d_{X/k} \in \text{III}(k, J)$, then $\langle \mathrm{Pic}^d_{X/k}, \mathrm{Pic}^d_{X/k} \rangle = 0$. In particular, $J$ over $k$ is even.*

*Proof.* Let $g$ be the genus of $X$. If $g$ is odd, then part (b) of Theorem 2 of [JL1] implies that for any completion $k_v$, $\mathrm{Pic}^d(X_{k_v}) = \emptyset$ if and only if $\mathrm{Pic}^d_{X/k}(k_v) = \emptyset$. Thus if $g$ is odd and $\mathrm{Pic}^d_{X/k} \in \text{III}(k, J)$, then the integer $N$ of Theorem 11 is zero; so $\langle \mathrm{Pic}^d_{X/k}, \mathrm{Pic}^d_{X/k} \rangle = 0$.

From now on we assume $g$ is even. Temporarily let us suppose $k = \mathbf{Q}$ and $d = 1$. By Corollary 4 of [JL1], $\mathrm{Pic}^1_{X/k} \in \text{III}(\mathbf{Q}, J)$. Since $g$ is even, $\mathrm{Disc}\, B = 2\ell$ where $\ell$ is a prime congruent to 3 or 5 modulo 8, as mentioned in the introduction of [JL1]. If $\ell \equiv 3 \pmod 8$, then [JL1] shows that $\mathrm{Pic}^1(X_{\mathbf{Q}_p}) = \emptyset$ exactly for $p = \ell$ and $p = \infty$. If $\ell \equiv 5 \pmod 8$, then [JL1] shows that $\mathrm{Pic}^1(X_{\mathbf{Q}_p}) = \emptyset$ exactly for $p = 2$ and $p = \infty$. In either case, $\langle \mathrm{Pic}^1_{X/k}, \mathrm{Pic}^1_{X/k} \rangle = 0$ by Theorem 11. Replacing $\mathbf{Q}$ by $k$ and $\mathrm{Pic}^1$ by $\mathrm{Pic}^d$ simply multiplies the result of the pairing by $d^2[k:\mathbf{Q}]$, so that the result is always zero. □

*Remarks.* 1. The results of [JL1] let one determine exactly when $\mathrm{Pic}^d_{X/k} \in \text{III}(k, J)$. We do not know, however, if $\mathrm{Pic}^d_{X/k}$ is trivial in every case in which it is everywhere locally trivial. Theorem 25 gives some evidence that this may be true.

2. Using our Corollary 12, Jordan and Livné [JL2] have recently shown that there are infinitely many Atkin-Lehner quotients of Shimura curves whose Jacobians over $\mathbf{Q}$ are odd.

3. The analogous questions for the standard modular curves are much easier. For $N \geq 3$ let $X(N)^{\mathrm{arith}}$ be the nonsingular projective model over $\mathbf{Q}$ of the affine curve $Y(N)^{\mathrm{arith}}$ that represents the functor "isomorphism classes of $\Gamma(N)^{\mathrm{arith}}$ elliptic curves" as on p. 482 of [Ka]. If $X$ is a curve over a number field $k$ that admits a map $X(N)^{\mathrm{arith}} \to X$ over $k$ for some $N \geq 3$, then $X(k)$ is nonempty, because there is the image of the cusp $\infty$ on $X(N)^{\mathrm{arith}}$, which is rational. Hence if $J$ is the Jacobian of $X$, the element $\mathrm{Pic}^n_{X/k} \in \text{III}(k, J)$ is trivial for every $n \in \mathbf{Z}$, and in particular $J$ is even.

### 10.2. Explicit examples.

PROPOSITION 26. *Let $g \geq 2$ and $t$ be integers with $g$ even. Let $X$ be the genus $g$ hyperelliptic curve*
$$y^2 = -(x^{2g+2} + x + t)$$
*over $\mathbf{Q}$ and let $J$ be its Jacobian. Then $J$ is odd if and only if $t > 0$.*

*Proof.* By Corollary 12, we need only count the deficient places of $X$. Since $g$ is even, and since $X$ has a degree 2 map to $\mathbf{P}^1$ over $\mathbf{Q}$, the number of such places is the same as the number of places $p$ for which $\mathrm{Pic}^1(X_{\mathbf{Q}_p}) = \emptyset$. We will show that $\mathrm{Pic}^1(X_{\mathbf{Q}_p}) \neq \emptyset$ for all finite $p$, and that the existence of a degree 1 line bundle for $p = \infty$ depends on the sign of $t$.

If $t$ is even, $x^{2g+2} + x + t$ has a zero in $\mathbf{Z}_2$, and if $t$ is odd, then $X$ has a $\mathbf{Q}_2$-rational point with $x$-coordinate $-1-t$, in each case by Hensel's lemma. Hence $X(\mathbf{Q}_2) \neq \emptyset$, so that $\mathrm{Pic}^1(X_{\mathbf{Q}_2}) \neq \emptyset$. If for an odd prime $p$ we have

$$-(x^{2g+2} + x + t) \equiv (a_{g+1}x^{g+1} + a_g x^g + \ldots + a_0)^2 \pmod{p}$$

with $a_{g+1}, a_g, \ldots, a_0 \in \overline{\mathbf{F}}_p$, then by equating high order coefficients we find $a_{g+1} \neq 0$, and then $a_i = 0$ for $i = g, g-1, g-2, \ldots, 1, 0$, contradicting the linear coefficient. Hence by Lemma 15, $\mathrm{Pic}^1(X_{\mathbf{Q}_p}) \neq \emptyset$.

If $t > 0$, then $t \geq 1$ and we see that $x^{2g+2} + x + t > 0$ for all $x \in \mathbf{R}$, by considering $x \leq -1$, $-1 \leq x \leq 0$, and $x \geq 0$ separately. Hence if $t > 0$, $X$ has no real divisor of odd degree. On the other hand, if $t \leq 0$, then $X(\mathbf{R})$ has a point with $x = 0$. □

PROPOSITION 27. *Let $J$ be the Jacobian of the genus 2 curve*

$$X : y^2 = -3(x^2+1)(x^2 - 6x + 1)(x^2 + 6x + 1)$$

*over $\mathbf{Q}$. Then $\mathrm{III}(J) \cong \mathbf{Z}/2\mathbf{Z}$.*

*Proof.* Let $E$ be the elliptic curve $64A1$ of [Cr] over $\mathbf{Q}$, which has CM by $\mathbf{Z}[i]$. Applying the results in [Ru] and using the data in Table 4 in [Cr], we find that $\mathrm{III}(E)(p)$ is trivial for all odd primes $p$. A 2-descent then shows $\mathrm{III}(E) = 0$. Another 2-descent shows that the 2-isogenous curve $64A3$ of [Cr],

$$E' : y^2 = (x-4)(x^2 + 4x - 28),$$

has $\mathrm{III}(E') = 0$.[18] Let $\psi$ be the $G_{\mathbf{Q}}$-module automorphism of $E'[2]$ that interchanges the nonrational 2-torsion points. Then the quotient of $E' \times E'$ by the graph of $\psi$ is a principally polarized abelian surface over $\mathbf{Q}$, and it equals $J$, by Proposition 3 in [HLP]. It follows that $\mathrm{III}(J) = \mathrm{III}(J)[2]$. Since $E$ has rank 0, so does $J$. A 2-descent on $J$ as in [PS] shows that its 2-Selmer group has $\mathbf{F}_2$-dimension 3, but $J(\mathbf{Q})/2J(\mathbf{Q})$ has $\mathbf{F}_2$-dimension only 2, so that $\mathrm{III}(J)[2] \cong \mathbf{Z}/2\mathbf{Z}$. See [St] for more on the implementation of the descent. □

*Remarks.* 1. The fact that $J$ is odd in Proposition 27 could have been predicted in advance by Corollary 12, since $X$ is deficient only at $p = 3$.

---

[18]The only reason for not applying Rubin's theorem directly to $E'$ is that the needed data in [Cr] are given for only one elliptic curve in each isogeny class, which happened to be $E$ in the case of interest.

2. Let $J$ be as in Proposition 27. One might ask whether the prediction given by the Birch and Swinnerton-Dyer conjecture for the value of $L(J,1)$ will be off by a power of 2, for instance. The answer is no, since the conjecture has been checked for $E$ and a theorem of Tate (see [Mi4, §1, Theorem 7.3]) implies that it holds for any abelian variety isogenous to $E \times E$.

PROPOSITION 28. *Let $J$ be the Jacobian of the genus* 2 *curve*

$$X : y^2 = -37(x^2+1)(5x^2-32)(32x^2-5)$$

*over* $\mathbf{Q}$. *Then* $\mathrm{III}(J)$ *is finite and* $\langle c,c \rangle = 0$, *but* $c \neq 0$. *Hence* $\mathrm{III}(J)$ *has square order, but* $\langle\ ,\ \rangle$ *is not alternating on it.*

*Proof.* We have that $J$ is $(2,2)$-isogenous over $\mathbf{Q}$ to $E \times E$ where $E$ is the rank 0 elliptic curve $y^2 = -37(x+1)(5x-32)(32x-5)$, isomorphic to curve $30A2$ in [Cr]. Kolyvagin's method [Ko] shows that $\mathrm{III}(E)$ is finite, and a 2-descent shows that $\mathrm{III}(E)[2] = 0$.[19] Hence $\mathrm{III}(J)$ is finite. Corollary 12 shows that $\langle c,c \rangle = 0$, since $X(\mathbf{Q}_p) \neq \emptyset$ for all $p$ including $\infty$.

To show $c \neq 0$, we perform a 2-descent on $J$ as in [PS]. Since $E$ has rank 0, $J$ has rank 0. We compute that the $(x-T)$ map injects $J[2]$ into $L^*/L^{*2}\mathbf{Q}^*$ (notation as in [PS]). Hence $(x-T) : J(\mathbf{Q})/2J(\mathbf{Q}) \to L^*/L^{*2}\mathbf{Q}^*$ is injective, and Theorem 11.3 of [PS] shows that $X$ cannot have a $\mathbf{Q}$-rational divisor of degree 1, i.e., that the homogeneous space $c = \mathrm{Pic}^1_{X/k}$ is nontrivial. The final statements follow from Corollaries 9 and 7. □

All our examples of odd Jacobians so far have been for curves of even genus. Next we give a family of genus 3 examples.

PROPOSITION 29. *Let $p$ be a prime with $p \equiv -1 \pmod{16}$. Let $X$ be the genus* 3 *curve over* $\mathbf{Q}$ *given in projective coordinates by*

$$X : x^4 + py^4 + p^2 z^4 = 0.$$

*Then the Jacobian of $X$ is odd.*

*Proof.* We use Corollary 12 and hope for an odd number of deficient places. Since $X$ has good reduction at all finite primes except possibly 2 and $p$, we need only check for deficiency at $\infty$, 2, and $p$. Any curve over $\mathbf{R}$ has real divisors of degree 2. Our curve $X$ also has a $\mathbf{Q}_2$-rational point that lifts from $(1:1:0)$ modulo $2^4$. It remains to prove that $X$ is deficient at $p$.

---

[19]The Birch and Swinnerton-Dyer conjecture predicts $\mathrm{III}(E) = 0$, and this could probably be proved unconditionally by making Kolyvagin's method explicit. This together with the 2-descent on $J$ would imply that $\mathrm{III}(J) \cong \mathbf{Z}/2\mathbf{Z} \times \mathbf{Z}/2\mathbf{Z}$.

CASSELS-TATE PAIRING   1141If not, then $X$ has a $\mathbf{Q}_p$-rational divisor of degree 2. Tripling the degree and applying Riemann-Roch yield an *effective* divisor of degree 6. Every partition of 6 includes a part that divides 6, so there exists an extension $L$ with $[L : \mathbf{Q}_p] = 6$ for which there exists a point $(x_0 : y_0 : z_0) \in X(L)$. Let $e$ and $f$ denote the ramification index and residue degree of $L/\mathbf{Q}_p$. If $e$ were odd, then $x_0^4$, $py_0^4$, $p^2 z_0^4$ would have distinct valuations modulo 4. Therefore $e = 2$ or $e = 6$. The valuation of $py_0^4$ is distinct from the others modulo 4, so $x_0^4$ and $p^2 z_0^4$ have the same valuation, and the valuation of $py_0^4$ is larger. Dividing by $x_0^4$, we find that $1 + p^2 z_0^4/x_0^4 \equiv 0$ modulo the maximal ideal; in particular $-1$ is a square in $\mathbf{F}_{p^f}$. This is impossible, since $p \equiv -1 \pmod 4$, and $f = 6/e$ is odd. □

*Remark.* Let $J$ be the Jacobian of $X$ in Proposition 29. We give here an unconditional proof that $\text{III}(J)$ is finite, for $p = 31$, 47, and 79.[20] First, $J$ is isogenous over $\mathbf{Q}$ to the product of the three CM elliptic curves

$$\begin{aligned} E : y^2 &= x^3 + x, \\ E' : y^2 &= x^3 + px, \\ E'' : y^2 &= x^3 + p^3 x. \end{aligned}$$

Rubin's results [Ru] let one prove that $\text{III}(E)$ is trivial. A numerical computation for $p = 31$, 47, and 79 shows that $\text{ord}_{s=1} L(E', s) = \text{ord}_{s=1} L(E'', s) = 1$. Then Kolyvagin's method [Ko] gives a proof that $\text{III}(E')$ and $\text{III}(E'')$ are finite. Thus $\text{III}(J)$ is definitely finite (and of nonsquare order) for these three $p$.

Finally we give some function field examples.

PROPOSITION 30. *Let $J$ be the Jacobian of the genus* 2 *curve*

$$X : y^2 = Tx^6 + x - aT$$

*over $\mathbf{F}_q(T)$, where $q$ is odd, and $a \in \mathbf{F}_q^* \setminus \mathbf{F}_q^{*2}$. If $\text{III}(J)$ is finite, then its order is twice a square.*

*Proof.* One easily checks that $X$ is deficient only at the infinite prime of $\mathbf{F}_q(T)$. □

*Remark.* There exists a (smooth, proper, geometrically integral) surface $Y$ over $\mathbf{F}_q$ equipped with a morphism $Y \to \mathbf{P}^1$ over $\mathbf{F}_q$ whose generic fibre is $X$. Since $Y$ is a rational surface, $\operatorname{Br}(Y)$ is finite [Mi1]. It is likely that the finiteness of $\operatorname{Br}(Y)$ can be proven to be equivalent to the finiteness of $\text{III}(J)$. The literature (see [Mi2] for instance) appears to contain a proof of the equivalence only under additional hypotheses which do not hold here, but it seems that the proof would extend without much difficulty to the general case.

---

[20]The same argument should apply to any given $p \equiv -1 \pmod{16}$; it relies only on the fact that the analytic ranks of $E'$ and $E''$ in what follows are less than or equal to 1. This can be checked for any given such $p$, and should hold for every such $p$. (For instance, it would follow from the weak form of the Birch and Swinnerton-Dyer conjecture, since a 2-descent shows that the algebraic ranks are less than or equal to 1.)

## 11. An open question of Tate about Brauer groups

Let $X$ be a (geometrically integral, smooth, projective) surface over a finite field $k$ of characteristic $p$. Because $X$ is smooth, the Brauer group $\operatorname{Br}(X)$ may be identified with $H^2_{\text{ét}}(X, \mathbf{G}_m)$.

Tate [Ta3, Th. 5.1] proved that for all primes $\ell \neq p$, $\operatorname{Br}(X)_{\text{nd}}(\ell)$ is finite, and admits a nondegenerate antisymmetric pairing with values in $\mathbf{Q}/\mathbf{Z}$. He asked whether this pairing was alternating for all $\ell \neq p$. Clearly the only case of interest is the case $\ell = 2 \neq p$. Tate's question remains unanswered in general, although Urabe [Ur] has made substantial progress: he has proved that the pairing is alternating on the subgroup $\ker \left( \operatorname{Br}(X)_{\text{nd}}(2) \to \operatorname{Br}(X_{\overline{k}})_{\text{nd}}(2) \right)$, from which he deduces that the answer to Tate's question is yes whenever $H^0(G_k, H^3_{\text{ét}}(X_{\overline{k}}, \mathbf{Z}_2(1)))_{\text{tors}} = 0$. Moreover, he has proved in general that the order of $\operatorname{Br}(X)_{\text{nd}}(2)$ is a square.

The same formalism as in our introduction and in Section 6 therefore gives an element $c \in \operatorname{Br}(X)_{\text{nd}}[2]$ with the property that $\langle x, x + c \rangle = 0$ for all $x \in \operatorname{Br}(X)_{\text{nd}}(2)$. Moreover Urabe's result about square order implies that $\langle c, c \rangle = 0$.

If the answer to Tate's question is yes, then $c$ should always be zero. Otherwise, one can hope for a more direct description of this canonically defined element of $\operatorname{Br}(X)_{\text{nd}}[2]$, analogous to our description of $c$ in the Shafarevich-Tate group case in Section 4, and perhaps this would facilitate the construction of a counterexample showing that the Brauer group pairing is not always alternating.

## 12. Appendix: Other definitions of the Cassels-Tate pairing

Here we present two more definitions of the Cassels-Tate pairing, and then prove the compatibility of all the definitions we have introduced. The "Albanese-Albanese definition" is new, but the "Weil pairing definition" is well-known (see [Mi4, p. 97]).

In the following, we always tacitly assume that all choices are made in such a way that all the pairings are defined. It is not difficult to see that this is always possible.

12.1. *The Albanese-Albanese definition.* Let $V$ and $W$ be varieties over a global field $k$ (conventions as usual), and suppose $D \in \mathrm{Div}(V \times W)$. One can simply mimic the construction in Section 3.2. Let $A = \mathrm{Alb}^0_{V/k}$ and $A' = \mathrm{Alb}^0_{W/k}$. We replace the second exact sequence in (3) by

$$0 \longrightarrow \mathcal{Y}^0(W_{k^{\mathrm{sep}}}) \longrightarrow \mathcal{Z}^0(W_{k^{\mathrm{sep}}}) \longrightarrow A'(k^{\mathrm{sep}}) \longrightarrow 0$$

and use the pairings

$$\mathcal{Y}^0(V_{k^{\mathrm{sep}}}) \times \mathcal{Z}^0(W_{k^{\mathrm{sep}}}) \longrightarrow k^{\mathrm{sep}*},$$
$$\mathcal{Z}^0(V_{k^{\mathrm{sep}}}) \times \mathcal{Y}^0(W_{k^{\mathrm{sep}}}) \longrightarrow k^{\mathrm{sep}*}$$

that agree on $\mathcal{Y}^0(V_{k^{\mathrm{sep}}}) \times \mathcal{Y}^0(W_{k^{\mathrm{sep}}})$ by Lang reciprocity. The same formalism as in Section 3.2 gives us a pairing

$$\langle \,,\, \rangle_D : \mathrm{III}(\mathrm{Alb}^0_{V/k}) \times \mathrm{III}(\mathrm{Alb}^0_{W/k}) \longrightarrow \mathbf{Q}/\mathbf{Z}.$$

The compatibility results we prove below will imply that $\langle \,,\, \rangle_D$ depends in fact only on the correspondence class of $D$, i.e., on the image of $D$ in $\mathrm{Pic}(V \times W)/(\pi_1^* \mathrm{Pic} V \oplus \pi_2^* \mathrm{Pic} W)$, where $\pi_1 : V \times W \to V$ and $\pi_2 : V \times W \to W$ are the projections.

12.2. *The Weil pairing definition.* Let $k$ be a global field of characteristic $p$ (we may have $p = 0$), and let $A$ be an abelian variety over $k$. Fix a positive integer $m$ with $p \nmid m$. Then the Weil pairing

$$e_m : A[m] \times A^\vee[m] \longrightarrow k^{\mathrm{sep}*}$$

can be defined as follows [La1, p. 173].[21] Given $a \in A[m]$ and $a' \in A^\vee[m]$, choose $\mathfrak{a} \in \mathcal{Z}^0(A_{k^{\mathrm{sep}}})$ and $\mathfrak{a}' \in \mathcal{Z}^0(A^\vee{}_{k^{\mathrm{sep}}})$ summing to $a$ and $a'$, respectively, and let[22]

$$e_m(a, a') := (m\mathfrak{a}) \cup \mathfrak{a}' - \mathfrak{a} \cup (m\mathfrak{a}'),$$

where $\cup$ is as in Section 12.1 for a Poincaré divisor $D \in \mathrm{Div}(A \times A^\vee)$. Let $\vee$ denote the cup-product pairing associated to $e_m$. We now use $e_m$ to define the Cassels-Tate pairing $\langle a, a' \rangle_{\mathrm{Weil}}$ for $a \in \mathrm{III}(A)$ and $a' \in \mathrm{III}(A^\vee)$ of order dividing $m$.

Choose $t \in H^1(A[m])$ and $t' \in H^1(A^\vee[m])$ mapping to $a$ and $a'$ respectively. Choose $\tau \in Z^1(A[m])$ and $\tau' \in Z^1(A^\vee[m])$ representing $t$ and $t'$ respectively. Choose $\sigma \in C^1(A[m^2])$ such that $m\sigma = \tau$ (i.e., so that $\sigma$ maps to $\tau$ under $A[m^2] \overset{m}{\to} A[m]$). Then $d\sigma$ takes values in $A[m]$, and $d\sigma \vee \tau'$ represents an element of $H^3(k^{\text{sep}*}) = 0$, so that $d\sigma \vee \tau' = d\bar{\varepsilon}$ for some $\bar{\varepsilon} \in C^2(k^{\text{sep}*})$.

---

[21]Our definition is the inverse of that in [La1]. On the other hand, our convention for the polarization associated to a line bundle is also opposite from that in [La1], so that we obtain the same answer to the question, "If $E$ is an elliptic curve over $\mathbf{C}$ corresponding to the lattice $\Lambda = \mathbf{Z} + \mathbf{Z}\tau$, if $E$ is identified with its dual, and if $P, Q \in E(\mathbf{C})$ correspond to $1/m, \tau/m \in \mathbf{C}/\Lambda$, then what is $e_m(P, Q)$?" That answer (unfortunately?) is $e^{-2\pi i/m}$.

[22]As usual, we write the group law in $k^{\mathrm{sep}*}$ additively.

Since $a_v = 0$, we can pick $\beta_v \in A(k_v^{\text{sep}})$ such that $d\beta_v$ equals the image of $\tau_v$ in $Z^1(G_v, A(k_v^{\text{sep}}))$. Choose $Q_v \in A(k_v^{\text{sep}})$ such that $mQ_v = \beta_v$. Let $\rho_v = dQ_v$ be reconsidered as an element of $Z^1(G_v, A[m^2])$. Note that $\sigma_v - \rho_v$ takes values in $A[m]$. Then

$$\bar{\gamma}_v := (\sigma_v - \rho_v) \vee \tau'_v - \bar{\varepsilon}_v \qquad \in C^2(G_v, k_v^{\text{sep}*})$$

is a 2-cocycle representing some

$$\bar{c}_v \in H^2(G_v, k_v^{\text{sep}*}) = \text{Br}(k_v) \xrightarrow{\text{inv}_v} \mathbf{Q}/\mathbf{Z}.$$

Define

$$\langle a, a' \rangle_{\text{Weil}} = \sum_{v \in M_k} \text{inv}_v(\bar{c}_v).$$

One checks that the value is well-defined, and unchanged if $m$ is replaced by a prime-to-$p$ multiple.

One advantage of this definition is that it extends to a motivic setting: see [Fl], which also outlines a proof that the definition is independent of choices made. A disadvantage is that the construction does not let one immediately define the pairing on the $p$-part of Ш, when $k$ has characteristic $p > 0$.

12.3. *Compatibility.* Let $V$ be a variety over a global field $k$, and let $A = \text{Alb}^0_{V/k}$. A choice of $P_0 \in V(k^{\text{sep}})$ gives rise to a morphism $\phi : V_{k^{\text{sep}}} \to A_{k^{\text{sep}}}$ that induces isomorphisms independent of $P_0$ that descend to $k$: $\phi_* : \text{Alb}^0_{V/k} \to \text{Alb}^0_{A/k}$ and $\phi^* : \text{Pic}^0_{A/k,\text{red}} \to \text{Pic}^0_{V/k,\text{red}}$. Our first result proves the equivalence of the Albanese-Picard pairings for $V$ and for $A$.

PROPOSITION 31. *For all $a \in $ Ш$(\text{Alb}^0_{V/k})$ and $a' \in $ Ш$(\text{Pic}^0_{A/k,\text{red}})$,*

$$\langle a, \phi^* a' \rangle_V = \langle \phi_* a, a' \rangle_A.$$

*Proof.* Let $\text{AlbMor}(V, A)$ denote the free group on the collection of Albanese morphisms $\psi : V_{k^{\text{sep}}} \to A_{k^{\text{sep}}}$ arising from choosing various basepoints $P_0$. There are then two ways of defining a $G_k$-equivariant trilinear pairing

$$\mathcal{Y}^0(V_{k^{\text{sep}}}) \times \text{AlbMor}(V, A) \times \text{Div}^0(A_{k^{\text{sep}}}) \longrightarrow k^{\text{sep}*},$$

either by using a generator $\psi \in \text{AlbMor}(V, A)$ to push forward 0-cycles from $V$ to $A$,

$$\mathcal{Y}^0(V_{k^{\text{sep}}}) \times \text{AlbMor}(V, A) \longrightarrow \mathcal{Y}^0(A_{k^{\text{sep}}}),$$

and then applying (2) for $A$, or by using $\psi$ to pull back divisors on $A$ to $V$, and then applying (2) for $V$. By the definition of (2), these trilinear pairings clearly agree. Similarly we have

$$\mathcal{Z}^0(V_{k^{\mathrm{sep}}}) \times \mathrm{AlbMor}(V, A) \times \frac{k^{\mathrm{sep}}(A)^*}{k^{\mathrm{sep}*}} \longrightarrow k^{\mathrm{sep}*},$$

and we may use $\cup$ for all the associated cup-product pairings without having to worry about ambiguity or nonassociativity.

Choose $\alpha \in Z^1(A(k^{\mathrm{sep}}))$ representing $a$, and lift $\alpha$ to $\mathfrak{a} \in C^1(\mathcal{Z}^0(V_{k^{\mathrm{sep}}}))$. Choose $\alpha' \in Z^1(A^\vee(k^{\mathrm{sep}}))$ representing $a'$, and lift $\alpha'$ to $\mathfrak{a}' \in C^1(\mathrm{Div}^0(A_{k^{\mathrm{sep}}}))$. We may then take $\phi \cup \mathfrak{a}'$ as the element of $C^1(\mathrm{Div}^0(V_{k^{\mathrm{sep}}}))$ required in the definition of $\langle a, \phi^* a' \rangle_V$, and $\mathfrak{a} \cup \phi$ as the element of $C^1(\mathcal{Z}^0(A_{k^{\mathrm{sep}}}))$ in the definition of $\langle \phi_* a, a' \rangle_A$. The $\eta$ in the definition of $\langle a, \phi^* a' \rangle_V$ then equals

$$d\mathfrak{a} \cup (\phi \cup \mathfrak{a}') - \mathfrak{a} \cup d(\phi \cup \mathfrak{a}')$$

and the corresponding element $\tilde{\eta}$ for $\langle \phi_* a, a' \rangle_A$ is

$$d(\mathfrak{a} \cup \phi) \cup \mathfrak{a}' - (\mathfrak{a} \cup \phi) \cup d\mathfrak{a}'.$$

A formal calculation shows that these are equal in $C^3(k^{\mathrm{sep}*})$, so we may take the $\varepsilon$'s in the two definitions to be the same.

Choose $\beta_v \in A(k_v^{\mathrm{sep}})$ such that $d\beta_v = \alpha_v$, and lift $\beta_v$ to $\mathfrak{b}_v \in \mathcal{Z}^0(V_{k_v^{\mathrm{sep}}})$, which we may push forward to $\mathfrak{b}_v \cup \phi \in \mathcal{Z}^0(A_{k_v^{\mathrm{sep}}})$ to serve as the corresponding "$\mathfrak{b}_v$" required in the definition of $\langle \phi_* a, a' \rangle_A$. The $\gamma_v$ in the definition of $\langle a, \phi^* a' \rangle_V$ then equals

$$(\mathfrak{a}_v - d\mathfrak{b}_v) \cup (\phi \cup \mathfrak{a}'_v) - \mathfrak{b}_v \cup d(\phi \cup \mathfrak{a}'_v) - \varepsilon_v$$

and the corresponding element for $\langle \phi_* a, a' \rangle_A$ is

$$((\mathfrak{a}_v \cup \phi) - d(\mathfrak{b}_v \cup \phi)) \cup \mathfrak{a}'_v - (\mathfrak{b}_v \cup \phi) \cup d\mathfrak{a}'_v - \varepsilon_v.$$

Again, these are formally equal. Summing invariants over $v$ completes the proof. □

Next we relate the Albanese-Picard definition to the homogeneous space definition.

PROPOSITION 32. *Let $A$ be an abelian variety over $k$, which we identify with $\mathrm{Alb}^0_{A/k}$. For all $a \in \mathrm{III}(A)$ and $a' \in \mathrm{III}(A^\vee)$,*

$$\langle a, a' \rangle_A = \langle a, a' \rangle.$$

*Proof.* Let $X$ be the homogeneous space of $A$ corresponding to $a$. Then $\mathrm{Alb}^0_{X/k} \cong A$, so that by Proposition 31, it will suffice to show

$$\langle a, a' \rangle_X = \langle a, a' \rangle.$$

Choose $P \in X(k^{\mathrm{sep}})$. For the definition of $\langle a, a' \rangle_X$, we may take $\mathfrak{a} = dP$, reconsidered as an element of $Z^1(\mathcal{Z}^0(X_{k^{\mathrm{sep}}}))$. Choose $\mathfrak{a}' \in C^1(\mathrm{Div}^0(X_{k^{\mathrm{sep}}}))$

representing $a'$. Then $d\mathfrak{a}'$ is the divisor of some $f' \in Z^2(k^{\text{sep}}(X)^*)$. The element $c_v \in H^2(k_v^{\text{sep}*})$ in the homogeneous space definition is obtained by evaluating $f'_v$ at any chosen $Q_v \in X(k_v)$. In the definition of $\langle a, a' \rangle_X$ on the other hand, we may take $\mathfrak{b}_v = P_v - Q_v \in \mathcal{Z}^0(X_{k_v^{\text{sep}}})$. Then $d\mathfrak{a} = d(dP) = 0$, so that

$$\eta = d\mathfrak{a} \cup \mathfrak{a}' - \mathfrak{a} \cup d\mathfrak{a}' = -dP \cup f' = d(-P \cup f')$$

and we may take $\varepsilon = -P \cup f'$. Next,

$$\begin{aligned}\gamma_v &= (\mathfrak{a}_v - d\mathfrak{b}_v) \cup \mathfrak{a}'_v - \mathfrak{b}_v \cup d\mathfrak{a}'_v - \varepsilon_v \\ &= (dP_v - d(P_v - Q_v)) \cup \mathfrak{a}'_v - (P_v - Q_v) \cup f'_v - (-P_v \cup f'_v) \\ &= Q_v \cup f'_v,\end{aligned}$$

and its cohomology class equals $c_v$, as desired. □

Now we relate the Albanese-Picard definition to the Albanese-Albanese definition.

PROPOSITION 33. *Let $V$ and $W$ be varieties over a global field $k$, and suppose $D \in \mathrm{Div}(V \times W)$. Then ${}^tD$ gives a homomorphism $\mathcal{Z}(W_{k^{\text{sep}}}) \to \mathrm{Div}(V_{k^{\text{sep}}})$, and we also denote by ${}^tD$ the induced homomorphism $\mathrm{Alb}^0_{W/k} \to \mathrm{Pic}^0_{V/k,\mathrm{red}}$. If $a \in \text{Ш}(\mathrm{Alb}^0_{V/k})$ and $a' \in \text{Ш}(\mathrm{Alb}^0_{W/k})$, then*

$$\langle a, a' \rangle_D = \langle a, {}^tDa' \rangle_V.$$

*Remark.* If $V = A$ is an abelian variety, if $W = A^\vee$, and if $D \in \mathrm{Div}(A \times A^\vee)$ is a Poincaré divisor, then the induced homomorphism ${}^tD : A^\vee \to A^\vee$ is the identity, and we find that

$$\langle a, a' \rangle_D = \langle a, a' \rangle_A = \langle a, a' \rangle.$$

*Proof.* This follows formally from the fact that ${}^tD$ maps the second exact sequence in the pair used to define $\langle a, a' \rangle_D$ down to the second exact sequence in the pair used to define $\langle a, {}^tDa' \rangle_V$:

$$\begin{array}{ccccccccc}0 & \longrightarrow & \mathcal{Y}^0(W_{k^{\text{sep}}}) & \longrightarrow & \mathcal{Z}^0(W_{k^{\text{sep}}}) & \longrightarrow & \mathrm{Alb}^0_{W/k}(k^{\text{sep}}) & \longrightarrow & 0 \\ & & \downarrow & & \downarrow & & \downarrow & & \\ 0 & \longrightarrow & \dfrac{k^{\text{sep}}(V)^*}{k^{\text{sep}*}} & \longrightarrow & \mathrm{Div}^0(V_{k^{\text{sep}}}) & \longrightarrow & \mathrm{Pic}^0(V_{k^{\text{sep}}}) & \longrightarrow & 0,\end{array}$$

the fact that the first exact sequences in the two pairs are the same, and the fact that these maps of exact sequences respect the two pairs of pairings needed in the definitions. □

Finally we relate the Weil pairing definition to the Albanese-Albanese pairing definition.

PROPOSITION 34. *Let $A$ be an abelian variety over a global field $k$, and let $\mathfrak{P} \in \mathrm{Div}(A \times A^{\vee})$ be a Poincaré divisor. If $a \in \mathrm{III}(A)$ and $a' \in \mathrm{III}(A^{\vee})$ have order prime to the characteristic of $k$, then*

$$\langle a, a' \rangle_{\mathfrak{P}} = \langle a, a' \rangle_{\mathrm{Weil}}.$$

*Proof.* We use notation consistent with that in Sections 12.1 and 12.2. Choose $t$, $t'$, $\tau$, $\tau'$, and $\sigma$ as in the latter. Choose $\mathfrak{s} \in C^1(G_k, \mathcal{Z}^0(A_{k^{\mathrm{sep}}}))$ representing $\sigma$. Then we may take $\alpha$ to be the image of $\tau$ in $Z^1(G_k, A(k^{\mathrm{sep}}))$, and take $\mathfrak{a} := m\mathfrak{s}$. Choose $\mathfrak{a}' \in C^1(G_k, \mathcal{Z}^0(A_{k^{\mathrm{sep}}}))$ representing $\tau'$. This determines $\eta$ and we choose $\varepsilon$ with $d\varepsilon = \eta$. Since

$$\begin{aligned}
d\sigma \vee \tau' - d\varepsilon &= [(md\mathfrak{s}) \cup \mathfrak{a}' - d\mathfrak{s} \cup (m\mathfrak{a}')] - [d\mathfrak{a} \cup \mathfrak{a}' - \mathfrak{a} \cup d\mathfrak{a}'] \\
&= [d\mathfrak{a} \cup \mathfrak{a}' - d\mathfrak{s} \cup (m\mathfrak{a}')] - [d\mathfrak{a} \cup \mathfrak{a}' - (m\mathfrak{s}) \cup d\mathfrak{a}'] \\
&= -d\mathfrak{s} \cup (m\mathfrak{a}') + \mathfrak{s} \cup d(m\mathfrak{a}') \\
&= -d(\mathfrak{s} \cup (m\mathfrak{a}')),
\end{aligned}$$

we may take $\bar{\varepsilon} = \varepsilon - \mathfrak{s} \cup (m\mathfrak{a}')$.

Now for the local choices. Choose $\beta_v$ such that $d\beta_v = \alpha_v$. (This is needed for both definitions.) Choose $Q_v$, which determines $\rho_v$. Choose $\mathfrak{q}_v \in \mathcal{Z}^0(A_{k_v^{\mathrm{sep}}})$ summing to $Q_v$. Then we may take $\mathfrak{b}_v = m\mathfrak{q}_v$. The difference between

$$\begin{aligned}
\bar{\gamma}_v &= (\sigma_v - \rho_v) \vee \tau'_v - \bar{\varepsilon}_v \\
&= (m(\mathfrak{s}_v - d\mathfrak{q}_v)) \cup \mathfrak{a}'_v - (\mathfrak{s}_v - d\mathfrak{q}_v) \cup (m\mathfrak{a}'_v) - \varepsilon_v + \mathfrak{s}_v \cup (m\mathfrak{a}'_v) \\
&= (\mathfrak{a}_v - d\mathfrak{b}_v) \cup \mathfrak{a}'_v + (d\mathfrak{q}_v) \cup (m\mathfrak{a}'_v) - \varepsilon_v
\end{aligned}$$

and

$$\begin{aligned}
\gamma_v &= (\mathfrak{a}_v - d\mathfrak{b}_v) \cup \mathfrak{a}'_v - \mathfrak{b}_v \cup d\mathfrak{a}'_v - \varepsilon_v \\
&= (\mathfrak{a}_v - d\mathfrak{b}_v) \cup \mathfrak{a}'_v - (m\mathfrak{q}_v) \cup d\mathfrak{a}'_v - \varepsilon_v
\end{aligned}$$

is $d(\mathfrak{q}_v \cup (m\mathfrak{a}'_v))$, which is a 2-coboundary, as desired. $\square$

## Acknowledgements

We thank J.-L. Colliot-Thélène, Torsten Ekedahl, Ehud Hrushovski, Bruce Jordan, Nick Katz, Hendrik Lenstra, Bill McCallum, Jim Milne, Peter Sarnak, Joe Silverman, and Jan van Geel for directing us to relevant references. We also thank Bill McCallum for sharing some of his unpublished notes elaborating the homogeneous space definition of the pairing, and J.-L. Colliot-Thélène for some comments and corrections.

UNIVERSITY OF CALIFORNIA, BERKELEY, CA, USA
*E-mail address*: poonen@math.berkeley.edu

MATHEMATISCHES INSTITUT DER HEINRICH-HEINE-UNIVERSITÄT, DÜSSELDORF, GERMANY
*E-mail address*: stoll@math.uni-duesseldorf.de